\documentclass{article}

\usepackage{amssymb,amsmath,theorem,euscript}

\newcounter{sec}

\newcounter{punct}[sec]

\def\punct{\refstepcounter{punct}{\arabic{sec}.\arabic{punct}.  }}

\def\COUNTERS{\addtocounter{sec}{1}
              \setcounter{punct}{0}
          \setcounter{equation}{0}
          \setcounter{theorem}{0}
                  }

\newtheorem{theorem}{Theorem}[sec]
\newtheorem{proposition}[theorem]{Proposition}

\newtheorem{lemma}[theorem]{Lemma}

\newtheorem{corollary}[theorem]{Corollary}

 \def\ov{\overline}
\def\wt{\widetilde}
 \newcommand{\rk}{\mathop {\mathrm {rk}}\nolimits}

\begin{document}

\def\OO{\mathrm{O}}
\def\GLO{\mathrm{GLO}}
\def\Coll{\mathrm{Coll}}
\def\kappa{\varkappa}
\def\Mat{\mathrm{Mat}}
\def\U{\mathrm U}

\def\R{\mathbb{R}}
\def\C{\mathbb{C}}

\def\la{\langle}
\def\ra{\rangle}

 \def\cA{\mathcal A}
\def\cB{\mathcal B}
\def\cC{\mathcal C}
\def\cD{\mathcal D}
\def\cE{\mathcal E}
\def\cF{\mathcal F}
\def\cG{\mathcal G}
\def\cH{\mathcal H}
\def\cJ{\mathcal J}
\def\cI{\mathcal I}
\def\cK{\mathcal K}
 \def\cL{\mathcal L}
\def\cM{\mathcal M}
\def\cN{\mathcal N}
 \def\cO{\mathcal O}
\def\cP{\mathcal P}
\def\cQ{\mathcal Q}
\def\cR{\mathcal R}
\def\cS{\mathcal S}
\def\cT{\mathcal T}
\def\cU{\mathcal U}
\def\cV{\mathcal V}
 \def\cW{\mathcal W}
\def\cX{\mathcal X}
 \def\cY{\mathcal Y}
 \def\cZ{\mathcal Z}
%%% END MATHCAL %%%%%%%%%%%%%%%%%%%%%%%%%%%%%%%%% %%%%%%%%%%%%%%%%%%%%%%%%%%%%%%%% %%%
\def\0{{\ov 0}}
% \def\1{{\ov 1}}
 %%%%%%%%%%%%%%%%%%%%%%%%%%%% %%%%%%%%%%%%%%%%%%%%%%%%%%%%%%%%%%% %%% BEGIN GOTIC
 \def\frA{\mathfrak A}
 \def\frB{\mathfrak B}
\def\frC{\mathfrak C}
\def\frD{\mathfrak D}
\def\frE{\mathfrak E}
\def\frF{\mathfrak F}
\def\frG{\mathfrak G}
\def\frH{\mathfrak H}
\def\frI{\mathfrak I}
 \def\frJ{\mathfrak J}
 \def\frK{\mathfrak K}
 \def\frL{\mathfrak L}
\def\frM{\mathfrak M}
 \def\frN{\mathfrak N} \def\frO{\mathfrak O} \def\frP{\mathfrak P} \def\frQ{\mathfrak Q} \def\frR{\mathfrak R}
 \def\frS{\mathfrak S} \def\frT{\mathfrak T} \def\frU{\mathfrak U} \def\frV{\mathfrak V} \def\frW{\mathfrak W}
 \def\frX{\mathfrak X} \def\frY{\mathfrak Y} \def\frZ{\mathfrak Z} \def\fra{\mathfrak a} \def\frb{\mathfrak b}
 \def\frc{\mathfrak c} \def\frd{\mathfrak d} \def\fre{\mathfrak e} \def\frf{\mathfrak f} \def\frg{\mathfrak g}
 \def\frh{\mathfrak h} \def\fri{\mathfrak i} \def\frj{\mathfrak j} \def\frk{\mathfrak k} \def\frl{\mathfrak l}
 \def\frm{\mathfrak m} \def\frn{\mathfrak n} \def\fro{\mathfrak o} \def\frp{\mathfrak p} \def\frq{\mathfrak q}
 \def\frr{\mathfrak r} \def\frs{\mathfrak s} \def\frt{\mathfrak t} \def\fru{\mathfrak u} \def\frv{\mathfrak v}
 \def\frw{\mathfrak w} \def\frx{\mathfrak x} \def\fry{\mathfrak y} \def\frz{\mathfrak z} \def\frsp{\mathfrak{sp}}
 %% This is Lie algebra %%% END GOTIC
%%%%%%%%%%%%%%%%%%%%%%%%%%%%%%%% %%%%%%%%%%%%%%%%%%%%%%%%%%%%%%%%%
%%% BEGIN MATHBF
 \def\bfa{\mathbf a} \def\bfb{\mathbf b} \def\bfc{\mathbf c} \def\bfd{\mathbf d} \def\bfe{\mathbf e} \def\bff{\mathbf f}
 \def\bfg{\mathbf g} \def\bfh{\mathbf h} \def\bfi{\mathbf i} \def\bfj{\mathbf j} \def\bfk{\mathbf k} \def\bfl{\mathbf l}
 \def\bfm{\mathbf m} \def\bfn{\mathbf n} \def\bfo{\mathbf o} \def\bfp{\mathbf p} \def\bfq{\mathbf q} \def\bfr{\mathbf r}
 \def\bfs{\mathbf s} \def\bft{\mathbf t} \def\bfu{\mathbf u} \def\bfv{\mathbf v} \def\bfw{\mathbf w} \def\bfx{\mathbf x}
 \def\bfy{\mathbf y} \def\bfz{\mathbf z} \def\bfA{\mathbf A} \def\bfB{\mathbf B} \def\bfC{\mathbf C} \def\bfD{\mathbf D}
 \def\bfE{\mathbf E} \def\bfF{\mathbf F} \def\bfG{\mathbf G} \def\bfH{\mathbf H} \def\bfI{\mathbf I} \def\bfJ{\mathbf J}
 \def\bfK{\mathbf K} \def\bfL{\mathbf L} \def\bfM{\mathbf M} \def\bfN{\mathbf N} \def\bfO{\mathbf O} \def\bfP{\mathbf P}
 \def\bfQ{\mathbf Q} \def\bfR{\mathbf R} \def\bfS{\mathbf S} \def\bfT{\mathbf T} \def\bfU{\mathbf U} \def\bfV{\mathbf V}
 \def\bfW{\mathbf W} \def\bfX{\mathbf X} \def\bfY{\mathbf Y} \def\bfZ{\mathbf Z} \def\bfw{\mathbf w}
 %%% END MATHBF
%%%%%%%%%%%%%%%%%%%%%%%%%%%%%%% %%%%%%%%%%%%%%%%%%%%%%%%%%%%%%%%%
 %%% BEGIN MATHBB
 \def\R {{\mathbb R }} \def\C {{\mathbb C }} \def\Z{{\mathbb Z}} \def\H{{\mathbb H}} \def\K{{\mathbb K}}
 \def\N{{\mathbb N}} \def\Q{{\mathbb Q}} \def\A{{\mathbb A}} \def\T{\mathbb T} \def\P{\mathbb P} \def\G{\mathbb G}
 \def\bbA{\mathbb A} \def\bbB{\mathbb B} \def\bbD{\mathbb D} \def\bbE{\mathbb E} \def\bbF{\mathbb F} \def\bbG{\mathbb G}
 \def\bbI{\mathbb I} \def\bbJ{\mathbb J} \def\bbL{\mathbb L} \def\bbM{\mathbb M} \def\bbN{\mathbb N} \def\bbO{\mathbb O}
 \def\bbP{\mathbb P} \def\bbQ{\mathbb Q} \def\bbS{\mathbb S} \def\bbT{\mathbb T} \def\bbU{\mathbb U} \def\bbV{\mathbb V}
 \def\bbW{\mathbb W} \def\bbX{\mathbb X} \def\bbY{\mathbb Y} \def\kappa{\varkappa} \def\epsilon{\varepsilon}
 \def\phi{\varphi} \def\le{\leqslant} \def\ge{\geqslant}

\def\GL{\mathrm {GL}}
\def\bGL{\mathbf {GL}}
\def\GLB{\mathrm {GLB}}

\def\bGr{\mathbf {Gr}}
\def\Gr{\mathrm {Gr}}
\def\bFl{\mathbf {Fl}}

\def\1{\mathbf {1}}

 \newcommand{\Dim}{\mathop {\mathrm {Dim}}\nolimits}
  \newcommand{\codim}{\mathop {\mathrm {codim}}\nolimits}
   \newcommand{\im}{\mathop {\mathrm {im}}\nolimits}
\newcommand{\ind}{\mathop {\mathrm {ind}}\nolimits}
\newcommand{\graph}{\mathop {\mathrm {graph}}\nolimits}

\def\F{\bbF}

\def\sm{\smallskip}

\begin{center}
\Large\bf
The space $L^2$ on semi-infinite Grassmannian over finite field

\large \sc
Yury A. Neretin%
\footnote{Supported by the grant FWF, P22122.}
\end{center}

{\small We construct a $\GL$-invariant measure on a semi-infinite Grassmannian 
over a finite field, describe the natural group of symmetries of this measure,
and decompose the space $L^2$ over the Grassmannian on irreducible representations.
The spectrum is discrete, spherical functions on the Grassmannian are given in terms
of the Al~Salam--Carlitz orthogonal polynomials. We also construct an invariant measure on the corresponding space
of flags.}

\section{Results of the paper}

\COUNTERS

{\bf\punct Group $\bGL(\ell\oplus\ell^\diamond)$.}
Let $\F_q$ be a finite field with $q$ elements.  Denote by $\ell$ the direct sum of countable
number of $\F_q$. Denote by $\ell^\diamond$ the direct product of countable
number of $\F_q$. Consider the linear space $\ell\oplus\ell^\diamond$ equipped with 
the natural topology.  Denote by $\bGL(\ell\oplus\ell^\diamond)$
the group of all continuous invertible linear operators
$\begin{pmatrix}a&b\\c&d\end{pmatrix}$
in $\ell\oplus\ell^\diamond$. A description of this group is given below in Section 2.

There is a natural homomorphism  $\theta:\bGL(\ell\oplus\ell^\diamond)\to\Z$
defined by 
$$
\theta \begin{pmatrix}a&b\\c&d\end{pmatrix}:=\dim(\ell/a\ell)-\dim(\ker a)
=-\bigl(\dim(\ell^\diamond/d\ell^\diamond)-\dim(\ker d)\bigr)
.
$$
Denote by $\bGL^0(\ell\oplus\ell^\diamond)$ the kernel of this homomorphism.

\sm

%%%%%%%%%%%%%%%%%%%%%%%%%%%%%%%%%%%%%%%%%%%%%%%%%%%%%%%%%%%%%%%%%

{\bf\punct Semi-infinite Grassmannian.}  For any infinite matrix $T$ over $\F_q$ consider
the subspace in $\ell\oplus \ell^\diamond$ consisting of vectors
$v\oplus vT\in \ell\oplus \ell^\diamond$. Denote by $\cM$ the set of all subspaces
in $\ell\oplus \ell^\diamond$
having such form.  We say that a subspace $M$ in $\ell\oplus \ell^\diamond$
is {\it semi-infinite} if there exists $L\in\cM$ such that
$$
\alpha(L):=\dim L/(L\cap M)<\infty, \qquad \beta(L):=\dim M/(L\cap M)<\infty
.
$$
We say that 
$$
\Dim(L):=\alpha(L)-\beta(L)
.
$$
 is the {\it relative dimension} of $L$.
Denote by $\bGr$ the set of all semi-infinite subspaces,
by $\bGr^\alpha$ the set of subspaces of relative dimension $\alpha$. 
  For any $g\in \bGL(\ell\oplus\ell^\diamond)$
 and semi-infinite $L$ we have
 $$
 \Dim (Lg)=\Dim L+\theta (g) 
 .$$
Therefore, $\bGL^0(\ell\oplus\ell^\diamond)$ acts on each $\bGr^\alpha$. This action is transitive.

The definition and basic properties of the semi-infinite Grassmannian are discussed in Section 3.

\sm

%%%%%%%%%%%%%%%%%%%%%%%%%%%%%%%%%%%%%%%%%%%%%%%%%%%%%%%%%%%

{\bf\punct The invariant measure on Grassmannian.}

\begin{theorem}
{\rm a)} There exists a unique %{\rm(}
up to a scalar factor
%{\rm)} 
finite $\bGL^0(\ell\oplus\ell^\diamond)$-invariant Borel
measure $\mu$ on $\bGr^0$.

\sm

{\rm b)} The restriction of the measure $\mu$ 
to the cell $\cM\simeq \F_q^{\infty\times\infty}$
is the product measure on the countable product of $\F_q$.

\sm

{\rm c)} If we normalize $\mu$ by the condition $\mu(\cM)=1$,
then the total measure of the Grassmannian is
$$
\mu(\bGr^0)=\prod_{j=1}^\infty (1-q^{-j})^{-1}
.$$
\end{theorem}

These statements are proved in Section 4.

\sm

%%%%%%%%%%%%%%%%%%%%%%%%%%%%%%%%%%%%%%%%%%%%%%%%%%%%%%5

{\bf\punct An intertwining operator.}
For $L\in \bGr^0$ consider the set $\Sigma_L$ of all semi-infinite subspaces $K\in\bGr^0$
such that
$$
\dim L/(L\cap K)=1\qquad \dim K/(L\cap K)=1
.$$

\begin{proposition}
There is a unique probabilistic measure $\nu_L$ on $\Sigma_L$ invariant with respect to the stabilizer
of $L$ in $\bGL^0(\ell\oplus\ell^\diamond)$.
\end{proposition}

The statement is proved in Subsection 5.1.

\begin{theorem}
{\rm a)}
The operator
$$
\Delta f(L)=\int_{\Sigma_L} f(K) d\nu_L(K)
$$
is a bounded self-adjoint $\bGL^0(\ell\oplus\ell^\diamond)$-intertwining operator
in $L^2(\bGr^0,\mu)$.

\sm

{\rm b)} The spectrum of $\Delta$ is the set $\{1, q^{-1}, q^{-2}, q^{-3},\dots\}$.

\sm

{\rm c)} The representations of $\bGL^0(\ell\oplus\ell^\diamond)$ in the eigenspaces
\begin{equation}
\Delta  f(L)=q^{-j} f(L)
\label{eq:Delta-lambda}
\end{equation}
are irreducible. % and pairwise nonequivalent.
\end{theorem}

Statement a) is proved in Subsection 5.1; b),c) in Subsection 5.5, 

\sm

%%%%%%%%%%%%%%%%%%%%%%%%%%%%%%%%%%%%%%%%%%%%%%%%%%%

{\bf\punct Invariant functions on the Grassmannian.}
Denote by $\bfP\subset \bGL^0(\ell\oplus \ell^\diamond)$ the stabilizer
of the subspace $0\oplus \ell^\diamond$.

\begin{proposition}
 Orbits  of $\bfP$ on $\bGr^0$ are enumerated by
$$
k:=\dim \bigl(L\cap (0\oplus\ell^\diamond)\bigr)
$$
{\rm ($k=0$, 1, 2,\dots)}.
Measures of orbits are
$$
\frac {q^{-k^2}}
{\prod_{j=1}^k (1-q^{-j})^2}
.
$$
\end{proposition}

For a proof see Corollary 4.14.

\begin{theorem}
{\rm a)} Any irreducible $\bGL^0(\ell\oplus\ell^\diamond)$-subrepresentation
in $L^2(\bGr^0)$ contains a unique $\bfP$-invariant function.

\sm

{\rm b)} The $\bfP$-invariant function in the eigenspace
$\Delta f= q^{-\alpha} f$ as a function of $k$ is given by the Al-Salam-Carlitz
polynomial
$$
V^{(1)}_\alpha(q^k,q^{-1})=
(-1)^\alpha q^{\alpha(\alpha-1)/2} {}_2\phi_0
\left[
\begin{matrix}
q^\alpha, q^k\\ -\! \!-
\end{matrix};q^{-1}; q^{-\alpha}\right]
.$$

\end{theorem}

This statement is proved in Subsection 5.5.

\sm

%%%%%%%%%%%%%%%%%%%%%%%%%%%%%%%%%%%%%%%%%%%%%%%

{\bf\punct Measures on flags.}
Let $[\alpha,\beta]$ be a segment of $\Z$.
We denote by $\bFl[\alpha,\beta]$ the space of all flags in
$\ell\oplus\ell^\diamond$ of the form
$$
L_\alpha\subset L_{\alpha+1}\subset \dots \subset L_\beta, \qquad 
L_j\in\bGr^j.
$$
We also allow infinite segments $(-\infty,\infty)$, $[\alpha,\infty)$,
$[-\infty,\beta]$.

\begin{proposition}
For any segment $I\subset \Z$ there is a unique $\bGL^0(\ell\oplus\ell^\diamond)$-invariant measure 
on $\bFl[I]$.
\end{proposition}

This is proved in Section 6. 

\sm

%%%%%%%%%%%%%%%%%%%%%%%%%%%%%%%%%%%%%%%%

{\bf\punct Representations of groups of  infinite-dimensional matrices over finite fields.}
In the remaining part of Introduction we discuss relation of our results with existing
theories, initial standpoints for the present paper,
  and some problems arising from our construction.

Now there is a relatively well-developed and relatively well-understood
representation theory of infinite-dimensional classical groups and
infinite symmetric groups. In a strange way, this is not so for groups of infinite matrices
over finite fields. There were several attempts to extend the widely used approaches 
to the groups $\GL(\infty,\F_q)$. We briefly discuss published works on this topic.

\smallskip

1) Let $G$ be a  group. Consider the set $\cK(G)$ of all positive definite central functions 
$\phi$ such that $\phi(e)=1$. According E.Thoma \cite{Tho}, a character of $G$ is an extreme point
of the compact set $\cK(G)$. Characters can be regarded as traces of type $\mathrm{II}_1$
 factor-representations 
of $G$ in the sense of Murray--von Neumann. On the other hand (\cite{Olsh-symm}),
 'characters' coincide with  functions on $G\times G$, spherical with respect to the diagonal
subgroup $\simeq G$.

 Thoma obtained a nice description of all characters
of  the group of finitely supported permutations of $\N$.
  In 1976 H.I.Skudlarek  \cite{Sku} tried to extend Thoma approach
to 
$$
\GL(\infty,\F_q)=\cup_{n=1}^\infty \GL(n,\F_q)$$
 His result was negative, 
 the group $\GL(\infty,\F_q)$ has no non-trivial characters. 
 
 \smallskip

2) G.I.Olshanski in series of papers (see \cite{OlshGB}, \cite{Olsh-symm}) proposed 
a $(G,K)$-pair approach 
to infinite-dimensional groups based on imitation of Harish-Chandra modules and on multiplications
of double cosets. In  \cite{Olsh-semi} (1991) he tried to extend this technology
to $\GL(\F_q,\infty)$.  More precisely, G.Olshanski
considered representations 'admissible' in the following sense. For
$j=1$, $2$, \dots consider the subgroup $G_j$ in $\GL(\infty,\F_q)$ consisting of matrices
$\begin{pmatrix}\1&0\\ 0& d \end{pmatrix}$, where the unit matrix $\1$ has size $j$.
Consider representations, which have fixed vectors with respect
to a sufficiently small subgroup $G_j$.
It appears that this class of representations is poor and only a small part of the picture
existing for real groups was imitated for finite fields. Apparently, now it is reasonable to return
to this approach keeping in mind extension of $(G,K)$-pairs in \cite{Ner-faa}.

\smallskip

3) A.M.Vershik an S.V.Kerov (see \cite{VK1}, \cite{VK2}, and the recent work 
\cite{GKV}) proposed a completion $\GLB(n,\F_q)$ of $\GL(n,\F_q)$
consisting of arbitrary invertible matrices over $\F_q$ having a only a finite 
number of elements under the diagonal. The subgroup of upper triangular
matrices is compact and the whole group is locally compact (but it is not a group
of type I). They obtained various 
constructions of representations  and quasiinvariant measures on flag varieties for the group
$\GLB(n,\F_q)$.

Our group $\bGL^0(\ell\oplus\ell^\diamond)$ is larger than $\GLB(\infty,\F_q)$, the natural topology on $\bGL^0(\ell\oplus\ell^\diamond)$ 
is weaker than on $\GLB(\infty,\F_q)$. Therefore representations of $\bGL^0(\ell\oplus\ell^\diamond)$
are also representations of $\GLB(n,\F_q)$. However, the present work has no intersection
with \cite{VK2}.

\sm

4)  Other group of works is \cite{GO1}--\cite{GO3}  (some our standpoints arise
from \cite{GO3}),
their main topic is the infinite symmetric group, but also there are discussed the space on 
decreasing flags in $\ell$,
$$
\ell=V\supset V_1\supset V_2\supset \dots,
$$
and measures on the flag space invariant with respect to the group of triangular matrices.

\sm

{\bf \punct  Inverse limits of homogeneous spaces.} In the present paper 
we propose a 'new' version 
of infinite-dimensional groups over finite fields. 
The initial purpose  was to obtain an analog of the following family of constructions. 

\sm

First (see D.Pickrell, \cite{Pick}) consider complex Grassmannians $\Gr_{2n}^n$ of $n$-di\-men\-si\-onal
subspaces in $\C^{2n}=\C^n\oplus \C^n$.
Consider the space $\Mat(n)$ of $n\times n$ complex matrices.
Any  $Z\in\Mat(n)$ determines an operator $\C^n\oplus 0\to 0\oplus \C^n$.
The graph of this operator determines a subspace in $\C^n\oplus\C^n$,
i.e., an element of $\Gr_{2n}^n$. Thus we get a bijection between $\Mat(n)$
and an open dense subset in $\Gr_{2n}^n$.
 The $\U(2n)$-invariant
probabilistic measure on $\Gr_{2n}^n$ is $\nu_n=\sigma_n\cdot\det(\1+ZZ^*)^{-2n}$,
where $\sigma_n$ is a normalizing scalar. For any matrix $Z\in\Mat(n)$
consider its left upper corner of size $(n-1)$. Thus we get a map
$\Mat(n)\to\Mat(n-1)$, and therefore a map $\Gr_{2n}^n\to\Gr_{2n-2}^{n-1}$.
 By the $\U(2n-2)$-invariance, the pushforward of $\nu_n$ is $\nu_{n-1}$.
 Thus we get a chain
$$
\dots\leftarrow\Gr_{2n-2}^{n-1}\leftarrow \Gr_{2n}^n \leftarrow\dots
$$
Its inverse limit $\lim\limits_{\leftarrow}\Gr_{2n}^n  $
 is equipped with a probabilistic measure $\nu$ and the group 
$\U(2\infty) =\lim\limits_\to \U(2n)$ acts on this space by measure preserving transformations. 

In fact, there is a natural one-parametric family ('Hua measures') on the inverse limit, we take
probabilistic measures of the form
$$
\nu_{n,s}=\sigma_{n}^{(s)}\cdot\det(\1+ZZ^*)^{-2n-s}
,$$
they are consistent with maps of Grassmannians, and we take an inverse limit as $n\to\infty$.

\sm

Several  constructions of such type are known;

\sm

1) projective limits of 
 classical groups $\U(n)$, $\OO(n)$, $\mathrm{Sp}(n)$ and symmetric spaces 
 (including Grassmannians over $\R$, $\C$, $\H$),
see  \cite{Ner-hua} 
 (also \cite{Ner-gauss}, Section 2.10), \cite {BO};

\sm

2) projective limits of
symmetric groups $S_n$, \cite{KOV};

\sm

3) projective limits of  $p$-adic Grassmannians, \cite{Ner-p}.

\sm

For two cases, $\lim\limits_\leftarrow \U(n)$ and $\lim\limits_\leftarrow S_n$,
now there is substantial harmonic analysis on the inverse limit, see \cite{KOV}, \cite{BO}.

\sm

The space $\bGr^0$ is the analog of such constructions for finite fields.
Notice that in our case the map of finite sets
$\Gr_{2n}^n\to \Gr_{2n-2}^{n-1}$ is defined only partially
and our construction is not precisely a  projective limit.
Also, the pre-limit objects in 1)-2) are analysis  
on groups and symmetric spaces. In our case, the pre-limit objects
are representations of $\GL(2N,\F_q)$ induced
from maximal parabolics
(the real analog are Stein representations \cite{Ste}, they do not survive at
the infinite-dimensional limit), the explicit decomposition
of such representations was obtained by 
A.B.Zhornitsky and A.V.Zelevinsky \cite{ZZ} and by D.Stanton,
\cite{Sta}. 
In this case functions invariant with respect to the maximal parabolic 
are $q$-Hahn polynomials (see, e.g., \cite{Koe}).
The Carlitz--Al-Salam II polynomials  are degenerations of $q$-Hahn
polynomials. 

We use limit  pass $n\to\infty$ for  proofs of Proposition 1.4 and Theorem 1.5.a,  however
it seems that it is difficult to prove Theorems 1.3 and 1.5.b in this way
(the obstacle is a mixture of different Grassmannians, which appear in  Subsection \ref{ss:pi}).

 The group $\bGL^0(\ell\oplus\ell^\diamond)$ arises as the group of symmetries of the measure $\mu$.
 % Notice that the construction of the Grassmannian $\bGr^0$ is canonical.
   Also, we get  actions of the same 
group $\bGL^0(\ell\oplus\ell^\diamond)$ on flag spaces. There arises a problem of decomposition of 
$L^2(\bFl[I])$. On the other hand, the measures on $\bFl[I]$ allow to apply the 
usual parabolic induction for construction of unitary representations of 
$\bGL^0(\ell\oplus\ell^\diamond)$. Therefore, it seems that the  Green theory (see, e.g., \cite{Mac})  of $\GL(n,\F_q)$
must survive in some form on the level of  the infinite-dimensional  group $\bGL^0(\ell\oplus\ell^\diamond)$.

Next,  the group $\bGL^0(\ell\oplus\ell^\diamond)$ is similar to infinite symmetric and infinite dimensional
classical real groups from the following point of view.

\sm

{\bf\punct The group $\bGL^0(\ell\oplus\ell^\diamond)$ as a $(G,K)$-pair.}
The group $\bGL^0(\ell\oplus\ell^\diamond)$ is a $(G,K)$-pair in the sense of \cite{Ner-book}, Section 6.5.
We formulate without proof  three simple statements.

\begin{proposition}
The pair $(\bGL^0, \bfP)$ is spherical, i.e., for any irreducible unitary representation
of $\bGL^0(\ell\oplus\ell^\diamond)$ the subspace of $\bfP$-fixed vectors has dimension $\le 1$.
\end{proposition}

Consider the subgroup $\bfP_k\subset \bfP\subset \bGL^0(\ell\oplus\ell^\diamond)$ consisting of block
$(k+\infty)+(k+\infty)$-matrices % $\begin{pmatrix}a&b\\c&d\end{pmatrix}$
 of the form
$$
\begin{pmatrix}
\bf \1&0&0&b_{12}
\\
a_{21}&a_{22}&b_{21}&b_{22}
\\
0&0& \1& d_{12}\\
0&0&0&d_{22}
\end{pmatrix}
$$
The double coset space $\bfP_k\setminus \bGL^0/\bfP_k$ admits a natural structure of a semigroup.
The multiplication is given in the following way. We represent each block
$a$, $b$, $c$, $d$ of $\bGL^0(\ell\oplus\ell^\diamond)$ as a block matrices of size $(k+\infty)$. 
Then we set
\begin{multline*}
\begin{pmatrix}
a_{11}&a_{12} &b_{11}&b_{12}\\
a_{21}&a_{22} &b_{21}&b_{22}\\
c_{11}&c_{12} &d_{11}&d_{12}\\
c_{21}&c_{22} &d_{21}&d_{22}
\end{pmatrix}
\circ 
\begin{pmatrix}
p_{11}&p_{12} &q_{11}&q_{12}\\
p_{21}&p_{22} &q_{21}&q_{22}\\
r_{11}&r_{12} &t_{11}&t_{12}\\
r_{21}&r_{22} &t_{21}&t_{22}
\end{pmatrix}:=
\\
:=
\begin{pmatrix}
a_{11}&a_{12}&0 &b_{11}&b_{12}&0\\
a_{21}&a_{22}&0 &b_{21}&b_{22}&0\\
0&0&\1&0&0&0\\
c_{11}&c_{12}&0 &d_{11}&d_{12}&0\\
c_{21}&c_{22}&0 &d_{21}&d_{22}&0\\
0&0&0&0&0&\1
\end{pmatrix}
\begin{pmatrix}
p_{11}&0&p_{12} &q_{11}&0&q_{12}\\
0&\1&0&0&0&0\\
p_{21}&0&p_{22} &q_{21}&0&q_{22}\\
r_{11}&0&r_{12} &t_{11}&0&t_{12}\\
0&0&0&0&\1&0\\
r_{21}&0&r_{22} &t_{21}&0&t_{22}
\end{pmatrix}
\end{multline*}
(in the right-hand side we have the usual product of matrices).

\begin{proposition}
 The $\circ$-multiplication is a well-defined associative operation
 on double cosets
$\bfP_k\setminus \bGL^0/\bfP_k$.
\end{proposition}

\begin{proposition}
Denote by $H_k$ the space of $\bfP_k$-fixed vectors in $L^2(\bGr^0)$. Denote by
$\Pi_k$ the orthogonal projection operator to $H_k$. Let
$\gamma\in \bfP_k\setminus \bGL^0/\bfP_k$ be a double coset. Let $g\in \bGL^0(\ell\oplus\ell^\diamond)$
be an element of $\gamma$. We define an operator $\pi(\gamma): H_k\to H_k$
by
$$
\pi(\gamma)f(L)=\Pi_k f(Lg)\bigr|_{H_k}
.
$$
Then for any double cosets $\gamma_1$, $\gamma_2\in \bfP_k\setminus \bGL^0/\bfP_k$
we have
$$
\pi(\gamma_1)\pi(\gamma_2)=\pi(\gamma_1\circ \gamma_2)
.
$$
\end{proposition}

There arises a natural conjecture: the last statement holds for arbitrary unitary representation of
the group $\bGL^0(\ell\oplus\ell^\diamond)$.

%It seems that the group $\bGL^0$ is a good infinite-dimensional analog of finite groups
%$\GL(n,\F_q)$. 

\sm

{\bf\punct Notation.}

--- $\F_p^\times$ is the multiplicative group
of the field $\F_q$.

--- $\#S$ is the number of elements of a finite set $S$.

%%%%%%%%%%%%%%%%%%%%%%%%%%%%%%%%%%%%%%%%%%%%%%%%%%%%%%%%%%%%%%%%%%%%%%
%%%%%%%%%%%%%%%%%%%%%%%%%%%%%%%%%%%%%%%%%%%%%%%%%%%%%%%%%%%%%%%%%%%
%%%%%%%%%%%%%%%%%%%%%%%%%%%%%%%%%%%%%%%%%%%%%%%%%%%%%%%%%%%%%%%%%%%%%%%
%%%%%%%%%%%%%%%%%%%%%%%%%%%%%%%%%%%%%%%%%%%%%%%%%%%%%%%%%%%%%%%%%%%%%%

\section{Linear operators in the  spaces $\ell$, $\ell^\diamond$, $\ell\oplus\ell^\diamond$}

\COUNTERS

Here we examine the group $\bGL(\ell\oplus\ell^\diamond)$ of invertible linear operators
in $\ell\oplus\ell^\diamond$.
 The main statement is
Theorem \ref{th:group1} about generators of $\bGL(\ell\oplus\ell^\diamond)$. 
We also get toy analogs
of Fredholm index theory and of  Banach's inverse operator theorem.

\hfill

{\bf\punct The spaces $\ell$ and $\ell^\diamond$.}
Denote by $\ell$ the direct sum of countable
number of $\F_q$. In other words, $\ell$ is the linear space of formal linear combinations
$$v=v_1 e_1+ v_2 e_2+\dots,$$
where $e_j$ are basis vectors
and $v_j\in \F_q$ are 0 for all but finite number of $v_j$.
 We equip
$\ell$ with the discrete topology.

Denote by $\ell^\diamond$ the direct product of countable
number of $\F_q$. In other words, $\ell$ is the linear space of all  formal linear combinations
$$w=w_1 f_1+ w_2 f_2+\dots,$$
where $f_j$ are basis elements and sequences $w_j\in\F_q$ are arbitrary.
 We equip $\F_q$ with the discrete topology. Then the direct product $\ell^\diamond$ becomes a compact
totally disconnected topological
space. 
We have obvious embedding $\ell\to \ell^\diamond$, evidently the image is dense.

The spaces $\ell$ and $\ell^\diamond$ are dual in the following sense:

\sm

--- any linear functional on $\ell$ has the form 
$h_w(v)=\sum v_j w_j$ for some $w\in \ell^\diamond$;

\sm

--- any continuous linear functional on $\ell^\diamond$ has the form 
$h_v(w)=\sum w_j v_j$ for some $w\in \ell^\diamond$.

\begin{lemma}
Any closed subspace in $\ell^\diamond$ is an annihilator of some subspace in $\ell$.
\end{lemma}

{\sc Proof.}
Here can refer to Pontryagin duality,  see, e.g., \cite{Mor}
Theorem 27.

\sm

%%%%%%%%%%%%%%%%%%%%%%%%%%%%%%%%%%%%%%%%%%%%%%%%%%%%%%%%%%%%%%

{\bf \punct Linear operators in $\ell$ and $\ell^\diamond$.}
Any  infinite matrix $A$ over $\F_q$ determines a linear operator
$\ell\to\ell^\diamond$ given by
$v\mapsto vA$,
\begin{multline*}
\begin{pmatrix} v_1&v_2&\dots\end{pmatrix}
\begin{pmatrix} a_{11}&a_{12}&\dots\\
a_{21}&a_{22}&\dots\\
\vdots&\vdots&\ddots
\end{pmatrix}
=\\=
\begin{pmatrix}
v_1 a_{11}+ v_2 a_{21}+\dots& v_1 a_{21}+ v_2 a_{22}+\dots&\dots
\end{pmatrix}
\end{multline*}
Here we write elements of $\ell$ as row-matrices.

\begin{lemma}
\label{l:ell-ell}
A matrix $A$ determines a linear operator
 $\ell\to\ell$ if and only if
 each row of $A$ contains only finite number of non-zero elements.
\end{lemma}

This is obvious (rows of $A$ are images of the basis vectors $e_j$). 
By duality we get the following statement.
 
\sm

\begin{lemma}
\label{l:ellcirc-ellcirc}
A matrix $A$ determines a continuous linear operator
 $\ell^\diamond\to\ell^\diamond$ if and only if 
 each column of $A$ contains only finite number of non-zero elements.
\end{lemma}

The transposition $A\mapsto A^t$ is an antiisomorphism ($(AB)^t=B^tA^t$) of the ring
of linear operators in $\ell$ and the ring of continuous linear operators in $\ell^\diamond$.

\begin{lemma}
\label{l:ellcirc-ell}
A matrix $A$ determines a linear operator
 $\ell^\diamond\to \ell$ if and only if
  $A$ contains only finite number of non-zero elements.
\end{lemma}

{\sc Proof.} Let $A$ be an operator $\ell^\diamond\to\ell$.
 By Lemma \ref{l:ellcirc-ellcirc} each column of $A$ contains only finite number of non nonzero
 elements. The transposed operator also acts $\ell^\diamond\to\ell$, therefore each row
 of $A$ contains only finite number of non-zero elements. Assume that $A$ contains 
 an infinite number of nonzero elements. Consider a sequence $f_{k_1}$, $f_{k_2}$,
 \dots such that $A f_{k_j}\ne 0$. Set of nonzero coordinates of each $A f_{k_j}$
 is finite, therefore we can choose a subsequence $A f_{k_{j_m}}$ such that the sets 
 of nonzero coordinates are mutually disjoint. Then
 $
 A \sum_m f_{k_{j_m}}\notin\ell
 $.
 \hfill $\square$
 
\sm

%%%%%%%%%%%%%%%%%%%%%%%%%%%%%%%%%%%%%%%%%%%%%%%%%%%%%%%%%%%%%%%%%%

{\bf\punct Groups of invertible operators.} Denote by $\bGL(\ell)$ and $\bGL(\ell^\diamond)$
groups of invertible continuous linear operators in $\ell$ and $\ell^\diamond$ respectively.

\begin{proposition}
\label{pr:ell-ell}
A linear operator $A:\ell\to\ell$  
{\rm(}see Lemma {\rm \ref{l:ell-ell}}{\rm)} is invertible if and only if it satisfies
the following additional conditions:

\sm

--- any finite collection of rows of the matrix $A$ is linear independent;

\sm

--- columns of $A$ are linear independent {\rm(}i.e., all finite and countable linear combinations
of columns are nonzero{\rm)}.
\end{proposition}

We emphasize that by Lemma \ref{l:ell-ell} {\it any countable linear combinations of columns
are well-defined}.

{\sc Proof.}
The first condition is equivalent to $\ker A=0$.
The second condition means that a linear functional on $\ell$ vanishing
on $\im A$ is zero, i.e. $\im A=\ell$. Under these two conditions
the operator $A$ is invertible.
\hfill $\square$

\sm

By the duality we get the following corollary

\begin{proposition}
\label{pr:ellcirc-ellcirc}
A linear operator $A:\ell^\diamond\to\ell^\diamond$  
{\rm(}see Lemma {\rm \ref{l:ellcirc-ellcirc}}{\rm)} is invertible if only if it satisfies
the following additional conditions:

\sm

--- any finite collection of columns of the matrix $A$ is linear independent

\sm

--- rows of $A$ are linear independent {\rm(}again, we admit  countable linear combinations
of rows{\rm)}.
\end{proposition}

%%%%%%%%%%%%%%%%%%%%%%%%%%%%%%%%%%%%%%%%%%%%%%%%%

{\bf\punct Fredholm operators in $\ell$.} We say that an operator
$A:\ell\to\ell$ is a {\it Fredholm operator} if 
$$
\dim\ker A<\infty \qquad  \codim \im A<\infty
.
$$
Its {\it index} is 
$$
\ind A:=\codim \im A-\dim\ker A
.
$$ 

\begin{lemma}
\label{l:Fcanonical}
An operator $A:\ell\to \ell$ is Fredholm if and only if
it can be represented in the form
$$
A=g_1 J g_2, 
$$
where $g_1$, $g_2\in\bGL(\ell)$, and $J$ is a
$(\alpha+\infty)\times(\beta+\infty)$ block matrix of the form
$$
\begin{pmatrix}0&0\\0&\1\end{pmatrix}
$$
with finite $\alpha$ and $\beta$.

Moreover, 
$$
\beta=\dim\ker A\qquad \alpha=\codim \im A
$$
\end{lemma}

{\sc Proof.} 'If' is obvious. 

Let $A$ be Fredholm. We take an arbitrary  operator $h_2\in\bGL(\ell)$ sending $\ker A$ to 
a subspace of the form $\oplus_{j<\beta} \F_q e_j$. Next, take $h_1\in\bGL(\ell^\diamond)$
sending the annihilator of $\im A$ to $\oplus_{j<\alpha} \F_q f_j$. Then
$ h_1^t A h_2^{-1}$ has the form $\begin{pmatrix}0&0\\0&B\end{pmatrix}$ with invertible
$B$.
\hfill $\square$

\begin{proposition}
Let $A:\ell\to\ell$ be Fredholm. Let $K:\ell\to\ell$
 be a continuous operator of finite rank. Then $A+K$
is a Fredholm operator of the same index.
\end{proposition}

{\sc Proof.} Without loss of generality we can assume that $A$ has canonical form
$\begin{pmatrix}0&0\\0&\1\end{pmatrix}$.
Since $\rk K$ is finite, it has only finite number of nonzero columns.
Therefore, without loss of generality, we can assume that $K$ has only one nonzero column.
After this a verification is straightforward.
\hfill $\square$

\begin{proposition}
\label{pr:Falmost}
Let $A$ be an operator $\ell\to\ell$. Assume that there are operators $B$, $C:\ell\to\ell$
such that
$$
AB=\1+K,\qquad CA=\1+L,
$$
where $K$ and $L$ have finite ranks. Then $A$ is Fredholm.
\end{proposition}

{\sc Proof.} Indeed,
$$
\ker A\subset \ker (\1+K),\qquad \im A\supset \im (\1+L).
$$

\begin{proposition}
\label{pr:Fdilatation}
An operator $A:\ell\to\ell$ is Fredholm if and only if there 
are finite $\alpha$, $\beta$ and an invertible matrix {\rm ('dilatation')}
$$
\begin{pmatrix}p&q\\r&A\end{pmatrix}
$$
of size $(\alpha+\infty)\times(\beta+\infty)$. Moreover,
$$
\ind(A)=\alpha-\beta.
$$
\end{proposition}

{\sc Proof.} Arrow $\Rightarrow$ easily follows from  Lemma \ref{l:Fcanonical}. Let us
prove $\Leftarrow$.

Represent $\begin{pmatrix}p&q\\r&A\end{pmatrix}^{-1}$ as a block
$(\beta+\infty)\times(\alpha+\infty)$ matrix
 $\begin{pmatrix}x&y\\z & B \end{pmatrix}$.
 Then
 $$
 AB=\1-ry,\qquad BA=\1-zq,
 $$
and we apply Proposition \ref{pr:Falmost}.
\hfill $\square$

\begin{proposition}
\label{pr:Fproduct}
Let $A$, $B$ be Fredholm operators in $\ell$. Then $AB$ is Fredholm and 
$$
\ind(AB)=\ind(A)+\ind(B)
.$$
\end{proposition}

{\sc Proof.} We include $A$ and $B$ to invertible matrices according
Proposition \ref{pr:Fdilatation},
$$
\wt A= \begin{pmatrix}p&q\\r&A\end{pmatrix},\qquad
\wt B= \begin{pmatrix}x&y\\z & B \end{pmatrix}
.
$$
Let their sizes be $(\alpha+\infty)\times (\beta+\infty)$
and $(\gamma+\infty)\times(\delta+\infty)$. If $\alpha=\delta$, then we write
$$
\wt A \wt B=\begin{pmatrix}p&q\\r&A\end{pmatrix}\begin{pmatrix}x&y\\z & B \end{pmatrix}
=\begin{pmatrix}*&*\\  *& AB+ry \end{pmatrix}.
$$
The product has size $(\gamma+\infty)\times (\beta+\infty)$. Therefore
$AB+rz$ is Fredholm of index $\gamma-\beta$. Therefore, $AB$ is Fredholm of the same index.

If $\alpha>\delta$, we change $\wt B$ to
$$
\wt B'=
\begin{pmatrix}
\1&0&0\\
0&x&y\\
0&z & B
 \end{pmatrix},
 $$
where the unit block has size $\alpha-\delta$ and multiply $\wt A\wt B'$.
If $ \alpha<\delta$, then we enlarge $\wt A$.
\hfill $\square$

\sm

%%%%%%%%%%%%%%%%%%%%%%%%%%%%%%%%%%%%%%%%%%%%%%%%%%%%%

{\bf\punct Fredholm operators in $\ell^\diamond$.}
We say that an operator $A:\ell^\diamond\to\ell^\diamond$ is {\it Fredholm}
if $\ker A$ is finite-dimensional, $\im A$  has finite codimension%
\footnote{Since the space $\ell^\diamond$ is compact, $\im A$ is closed for any
continuous $A$}.
The {\it index} is $(\codim \im A-\dim\ker A)$ as above.

\begin{proposition}
An operator $A:\ell^\diamond\to\ell^\diamond$ is Fredholm if and only if 
$A^t$ is Fredholm in $\ell$,
$$
\ind A^t=-\ind A
$$
\end{proposition}

%%%%%%%%%%%%%%%%%%%%%%%%%%%%%%%%%%%%%%%%%%%%%%%%%%%%%%%%%

{\bf\punct Linear operators in $\ell\oplus\ell^\diamond$.} Now consider the space
$\ell\oplus \ell^\diamond$, denote  the standard basis in $\ell$ by 
$e_j$, and the standard basis in $\ell^\diamond$ by $f_j$. 

The dual space is $\ell^\diamond\oplus \ell$. Note that $\ell\oplus\ell^\diamond$
is a locally compact Abelian group, therefore we have a possibility to refer
to the Pontryagin duality. In particular,

\begin{lemma}
\label{l:Pdual2}
Any closed subspace in $\ell\oplus\ell^\diamond$ is the annihilator of
a closed subspace in $\ell^\diamond\oplus\ell$.
\end{lemma}

We write continuous linear operators in $\ell\oplus\ell^\diamond$ as block 
$(\infty+\infty)\times (\infty+\infty)$
matrices  
$$g=\begin{pmatrix}a&b\\c&d\end{pmatrix}$$
 (and regard elements of
$\ell\oplus\ell^\diamond$ as block row matrices).
We also use the notation
$$
g=\begin{pmatrix}g_{11}&g_{12}\\g_{21}&g_{22}\end{pmatrix}
$$

By Lemmas \ref{l:ell-ell}--\ref{l:ellcirc-ell} such matrices satisfy the conditions:

\sm

a) Any row of $a$ contains finite number of nonzero elements;

\sm

b) Any column of $d$ contains finite number of nonzero elements;

\sm

c) $c$ has only finite number of nonzero elements.

\sm

{\bf\punct The group $\bGL(\ell\oplus\ell^\diamond)$.}
We denote by $\bGL(\ell\oplus\ell^\diamond)$
the group of all continuous invertible linear operators
in $\ell\oplus\ell^\diamond$.

\begin{lemma}
If $\begin{pmatrix}a&b\\c&d\end{pmatrix}\in\bGL(\ell\oplus\ell^\diamond)$, then $a$ and $d$
are Fredholm.
\end{lemma}

{\sc Proof.} Let $\begin{pmatrix}p&q\\r&t\end{pmatrix}$ be the inverse matrix. Then
\begin{align*}
\begin{pmatrix}\1&0\\0&\1\end{pmatrix}=
\begin{pmatrix}a&b\\c&d\end{pmatrix} \begin{pmatrix}p&q\\r&t\end{pmatrix}
=\begin{pmatrix}ap+br& *\\ *&cq+dt\end{pmatrix};
\\
\begin{pmatrix}\1&0\\0&\1\end{pmatrix}=
\begin{pmatrix}p&q\\r&t\end{pmatrix}
\begin{pmatrix}a&b\\c&d\end{pmatrix}=
\begin{pmatrix}pa+qc&*\\ *&rb+td\end{pmatrix}.
\end{align*} 
Matrices $d$, $r$ are finite, therefore $br$, $cq$, $qc$, $rb$ have finite rank.
By Proposition \ref{pr:Falmost}, $a$ and $d$ are Fredholm
\hfill $\square$
 
\sm 
 
{\bf\punct Generators of  $\bGL(\ell\oplus\ell^\diamond)$.}
Next, we define certain subgroups in $\bGL(\ell\oplus\ell^\diamond)$.

\sm

1) The group $\GL(2\infty,\F_q)$ is the group of finite invertible matrices, i.e., matrices
$g$ such that $g-1$ has only finite number of nonzero elements. We represent 
$\GL(2\infty,\F_q)$ as a union of subgroups $\GL(2n,\F_q)$
consisting of block $(n+\infty+n+\infty)$-matrices  of the form
$$
\begin{pmatrix}
u&0&v&0\\
0&\1&0&0\\
w&0&x&0\\
0&0&0&\1
\end{pmatrix}.
$$
The group $\GL(2\infty, \F_q)$ is the inductive limit
$$
\GL(2\infty, \F_q)=\lim_{n\to \infty} \GL(2n,\F_q)
.$$

2) The 'parabolic' subgroup $\bfP$ consists of matrices $\begin{pmatrix}a&b\\0&d\end{pmatrix}$
such that $a\in\bGL(\ell)$, $d\in\bGL(\ell^\diamond)$.

3) Next consider an operator $J$ defined by
\begin{align}
e_i J&=e_{i-1}, \text{where $i\ge 2$};
\nonumber
\\
e_1 J&=f_1;
\label{eq:J}
\\
f_k J&=  f_{k+1}.
\nonumber
\end{align}
Denote by $\cZ\subset \bGL(\ell\oplus\ell^\diamond)$ the cyclic
group generated by $J$.

\sm

Now we present two descriptions of the group $\bGL(\ell\oplus\ell^\diamond)$.

\begin{theorem}
\label{th:group1}
The group $\bGL(\ell\oplus\ell^\diamond)$ is generated by the subgroups $\GL(2\infty,\F_q)$,
$\bfP$, and $\cZ$.
\end{theorem}

{\sc Proof.} Let $g\in\bGL(\ell\oplus\ell^\diamond)$. We choose
$k$ such that for $g':=J^k g$ the block $g'_{11}$ has index 0. A multiplication
of the form
\begin{equation}
\begin{pmatrix}u&0\\0&\1\end{pmatrix}
\begin{pmatrix}g'_{11}&g'_{12}\\g'_{21}&g'_{22}\end{pmatrix}
\begin{pmatrix}v&0\\0&\1\end{pmatrix}:=g'',\qquad 
\text{where $u$, $v\in\bGL(\ell)$,}
\label{eq:step-1}
\end{equation}
allows to get
$$
g''_{11}=\begin{pmatrix}0&0\\0&\1 \end{pmatrix}
$$
(the size of the matrix is $(k+\infty)\times(k+\infty)$).
We write $g''_{21}$ as a block $\infty\times (k+\infty)$-matrix 
 $g''_{21}=\begin{pmatrix} p&q\end{pmatrix}$.
 We set
 $$
 g''':=\begin{pmatrix}\1&0&0\\ 0&\1&0\\0&-q&\1 \end{pmatrix} g''
 ,
 $$
the size of the first factor is $(k+\infty+\infty)\times (k+\infty+\infty)$, it
is contained in $\GL(2\infty,\F_q)$.
In this way, we  get
$$
g'''_{21}=\begin{pmatrix} p&0\end{pmatrix}
$$
(other blocks are the same). The rank of $p$ is $k$. We choose
$k$ linear independent rows of $p$ and permute them with first $k$ rows
of $g'''_{11}$ (this corresponds to a multiplication by a certain element
of $\GL(2\infty,\F_q)$). Thus we get a new matrix $g^{IV}$ with
$$
g^{IV}_{11}=\begin{pmatrix}h&0\\0&\1 \end{pmatrix}
$$
Repeating the operation (\ref{eq:step-1}), we can get 
$g^V_{11}=\begin{pmatrix}1&0\\0&\1\end{pmatrix}$,
 $g^V_{21}=\begin{pmatrix}p'&0\end{pmatrix}$.
After this we kill $p'$ by 
$$
g^{VI}:= \begin{pmatrix}\1&0&0\\ 0&\1&0\\-p'&0&\1 \end{pmatrix} g^V
$$
and get a matrix of the form
$$
g^{VI}=\begin{pmatrix}\1& x\\0&y \end{pmatrix}
.$$
This matrix is contained in $\bfP$. 
\hfill $\square$

\sm

{\bf\punct The inverse operator theorem.} The following statement does not used below.

\begin{theorem}
\label{th:group2}
 If a continuous linear operator $g$ in $\ell\oplus\ell^\diamond$ is bijective,
then $g$ is invertible.
\end{theorem}

{\sc Proof.} Since the block $c$ is finite, the image  of $0\oplus\ell^\diamond$
is contained in some subspace
$$
H=
\Bigl(\bigoplus_{j\le k} \F_q e_j\Bigr)\oplus \ell^\diamond
.
$$
The subspace $ \ell^\diamond g$ is compact, therefore it is a closed  in $H$.
For any closed subspace $L\subset H$ the dimension $H/L$ is finite  or continual.
On the other hand the dimension $(\ell\oplus \ell^\diamond) /\ell^\diamond$
is countable, therefore the dimension $(\ell\oplus \ell^\diamond)/ \ell^\diamond g$
is countable. Thus, $\dim (H/\ell^\diamond g)$ is finite. Applying
an appropriate element of $\GL(2\infty,\F_q)$ we can put $ \ell^\diamond g$ to position
$$
\Bigl(\bigoplus_{j\le m} \F_q e_j\Bigr)\oplus \ell^\diamond
\qquad \text{or}\quad \bigoplus_{i\ge \alpha}\F_q f_i
.
$$
Applying an appropriate power $J^s$ we can assume $\ell^\diamond g=\ell^\diamond$.
The map $g:\ell^\diamond\to\ell^\diamond g$ is continuous. Therefore the inverse
map $g^{-1}:\ell^\diamond g\to\ell^\diamond$ is continuous. Since
$(\ell\oplus\ell^\diamond)/\ell^\diamond g$ is discrete, the map $g^{-1}$ is automatically
continuous on the whole $\ell\oplus\ell^\diamond$.
\hfill $\square$

\begin{proposition}
\label{pr:group3}
Let $g$ be a continuous operator in $\ell\oplus\ell^\diamond$, $\ker g=0$, and
$\im g$ be dense. Assume that 
$$
\dim \ell^\diamond/(\ell^\diamond\cap \ell^\diamond g)
$$
is finite. Then $g$ is invertible.
\end{proposition}

{\sc Remark.} The last condition is necessary. Consider 
the spaces $V=\ell\oplus\ell\oplus \ell^\diamond$ and $W=\ell\oplus\ell^\diamond\oplus \ell^\diamond$
(both spaces can be identified with $\ell\oplus \ell^\diamond$). Consider 
the operator $A:V\to W$ defined by the matrix
$$
\begin{pmatrix}
0&0&\1\\
\1&0&0\\
0&\1&0
\end{pmatrix}
$$
Then $\ker A=0$ and $\im A:=\ell\oplus\ell^\diamond\oplus \ell$ is dense in $W$.
But $A$ is not invertible. %\hfill $\square$

\sm

{\sc Proof.} As in the proof of Theorem \ref{th:group2} we can assume 
$ \ell^\diamond g=\ell^\diamond$. Thus $(\ell\oplus\ell^\diamond)/\ell^\diamond g$ is 
a space with discrete topology. Since the image of $g$ is dense, $g$ must be surjective.
\hfill $\square$ 

\sm

{\bf\punct The subgroup $\bGL^0(\ell\oplus\ell^\diamond)$.}

\begin{theorem}
Let 
$\begin{pmatrix}a&b\\c&d \end{pmatrix}$ range in $\bGL(\ell\oplus\ell^\diamond)$.

\sm

{\rm a)} $\ind a+\ind d=0$.

\sm

{\rm b)} The map 
$$
\theta\begin{pmatrix}a&b\\c&d \end{pmatrix}=\ind a
$$ 
is a homomorphism $\bGL(\ell\oplus\ell^\diamond)\to \Z$.

\sm

{\rm c)}
The homomorphism $\theta$ is also determined by the conditions
\begin{align}
\theta(g)&=0, \qquad \text{if $g\in\bfP$};
\label{eq:theta1}
\\
\theta(h)&=0, \qquad \text{if $h\in\GL(2\infty,\F_q)$};
\\
\theta(J)&=1
\label{eq:theta3}
.\end{align}
\end{theorem}

{\sc Proof.} b) Let $\begin{pmatrix}p&q\\r&t\end{pmatrix}$
be the inverse matrix,
$$
\begin{pmatrix}a&b\\c&d\end{pmatrix} \begin{pmatrix}p&q\\r&t\end{pmatrix}
=\begin{pmatrix}ap+br& *\\ *& *\end{pmatrix}
.$$
By Proposition \ref{pr:Fproduct}, $ap$ is Fredholm and
$\ind ap=\ind a+\ind p$.
The matrix $r$ is finite, $br$ has finite rank. Therefore
$\ind (ap+br)=\ind ap$.

\sm

c) is obvious.

\sm

a) $\wt\theta(g):=-\ind d$ also is a homomorphism
$\bGL(\ell\oplus\ell^\diamond)\to \Z$. We have $\wt\theta(g)=\theta(g)$
on generators (\ref{eq:theta1})--(\ref{eq:theta3}). Therefore $\wt\theta=\theta$.
\hfill $\square$

\sm

We denote by $\bGL^0(\ell\oplus\ell^\diamond)$ the kernel
of the homomorphism $\theta$.

The proof of Theorem \ref{th:group1} gives also the following statement.

\begin{proposition}
\label{pr:cor}
Any element of $\bGL^0$ can be represented
  as
  $$
  g=h q  r, \qquad\text{where $h\in\bGL(\ell)$, $q\in\GL(2\infty,\F_q)$, $r\in\bfP$}.  
  $$
  \end{proposition}

%%%%%%%%%%%%%%%%%%%%%%%%%%%%%%%%%%%%%%%%%%%%%%%%%%%%%%%%%%%
%%%%%%%%%%%%%%%%%%%%%%%%%%%%%%%%%%%%%%%%%%%%%%%%%%%%%%%%%%%
%%%%%%%%%%%%%%%%%%%%%%%%%%%%%%%%%%%%%%%%%%%%%%%%%%%%%%%%%%%
%%%%%%%%%%%%%%%%%%%%%%%%%%%%%%%%%%%%%%%%%%%%%%%%%%%%%%%%%%%

\section{The semi-infinite Grassmannian}

\COUNTERS

{\bf\punct The Grassmannian.}
We use alternative notation for the following subspaces in $\ell\oplus\ell^\diamond$:
$$
V:=\ell\oplus 0, \qquad  W:=0\oplus \ell^\diamond
.
$$

 We define the {\it semi-infinite Grassmannian}
 $\bGr$ as the set of all subspaces $L$ in $\ell\oplus\ell^\diamond$
having the form
\begin{equation}
L=V g,\qquad 
\text{where $g=\begin{pmatrix}a&b\\c&d\end{pmatrix}\in \bGL(\ell\oplus\ell^\diamond)$}
\label{eq:def-gr}
,\end{equation}
i.e., $\bGr$ is the $\bGL(\ell\oplus\ell^\diamond)$-orbit of $\ell^\diamond$.

We define the {\it relative dimension} of $L$ as
$$
\Dim(L):=\ind(a)=\theta(g)
.$$

By definition,
$$
\Dim Lh=\Dim L+\theta(h), \qquad L\in\bGr, h\in\bGL(\ell\oplus\ell^\diamond).
$$

We denote by $\bGr^0$ the $\bGL^0(\ell\oplus\ell^\diamond)$-orbit of
$V$. It consists of elements of $\bGr$ having relative dimension 0. 

More generally, denote by $\bGr^j\subset \bGr$ the set of subspaces of relative dimension $j$

The following lemma implies that $\Dim L$ does not depend on the choice of $g$
in (\ref{eq:def-gr}).

\begin{lemma}
Let $L\in\bGr$. Denote by $p(L)$ the projection of $L$ to $V$ along $W$.
Then
$$
\Dim L= \dim (L\cap W)-\dim(V/p(L)).
$$
\end{lemma}

{\sc Proof.} The subspace  (\ref{eq:def-gr}) consists of vectors
\begin{equation}
va\oplus vb,\qquad \text{where $v$ ranges in $V$.}
\label{eq:v}
\end{equation}
We have 
$
L\cap W=(\ker a) b.
$
The operator $b$ does not vanish on $\ker a$, otherwise $\ker g\ne0$.
Therefore $\dim L\cap W=\dim \ker a$.

On the other hand $p(L)=\im a$.
\hfill $\square$

\sm

 Consider a linear operator $T:\ell\to\ell^\diamond$.
Consider its graph $\graph(T)\subset \ell\oplus\ell^\diamond$, i.e., the set of vectors
of the form $v\oplus vT$. Denote by $\cM$ the set of all such subspaces in $\ell\oplus\ell^\diamond$.

\begin{lemma}
\label{l:linear-fractional}
{\rm a)} $\cM\subset \bGr^0$.

\sm

{\rm b)} If  $a$ in {\rm(\ref{eq:def-gr})} is invertible, then 
$T=a^{-1} b$.

\sm

{\rm c)} If $a+Tc$ is invertible, then
\begin{equation}
\label{eq:linear-fractional}
\graph(T)\begin{pmatrix}a&b\\c&d\end{pmatrix}=
\graph\bigl( (a+Tc)^{-1}(b+Td)\bigr)
\end{equation}

{\rm d)} Let $a+Tc$ be invertible.
 Consider the images $\wt V$, $\wt W$ of subspaces $V$ and $W$ under 
$\begin{pmatrix}a&b\\c&d\end{pmatrix}$, i.e. spaces consisting of
vectors $va\oplus vc$ and $wb\oplus wc$, where $v$ ranges in $V$
and $w$ ranges in $W$. Then the subspace $\graph T\subset V\oplus W$
is a graph of an operator $\wt V\to \wt W$, whose matrix is $ (a+Tc)^{-1}(b+Td)$.
\end{lemma}

{\sc Proof.}
a) We set $g=\begin{pmatrix}\1&T\\0&\1\end{pmatrix}$.

\sm

b) In (\ref{eq:v}) we set $v'= av$ and get $v'\oplus v'a^{-1}b$.

\sm

c) We write
$$
\begin{pmatrix}v &vT \end{pmatrix} \begin{pmatrix}a&b\\c&d\end{pmatrix}
=\begin{pmatrix}v(a+Tc)&v(b+Td) \end{pmatrix}
$$
and set $v'= v(a+Tc)$.

\sm 

d) is a reformulation of c).

\hfill $\square$

\sm

%%%%%%%%%%%%%%%%%%%%%%%%%%%%%%%%%%%%%%%%%%%%%%%%%%%%%%%%%%%%%%%%%%%

{\bf\punct A characterization of elements of $\bGr$.}

\begin{lemma}
{\rm a)} Let $L=Vg\in\bGr$, $M$ be a subspace in $L$, $\dim L/M<\infty$.
Then $M\in\bGr$.

\sm

{\rm b)} Let $L\in\bGr$, let $N\supset L$ be  a subspace in $\ell\oplus \ell^\diamond$
such that $\dim N/L<\infty$.
Then $N\in \bGr$.

\sm

{\rm c)} If $L\supset M$ are elements of $\bGr$, then
$$
\Dim L=\Dim M+\dim L/M.
$$ 
\end{lemma}

{\sc Proof.} Let us prove a). $Mg^{-1}$ is a subspace in $V$. Consider an invertible operator $h$
in $V\simeq\ell$ such that 
$$
M  g^{-1}h^{-1} =\bigoplus_{j\ge k} \F_qe_j
$$
Then
$
M= V J^k h g 
$. 
\hfill $\square$

\begin{proposition}
\label{pr:LM}
For any subspace $L\in \bGr$ there exists $M\in \cM$ such that 
$\dim L/(L\cap M)<\infty$, $\dim M/(L\cap M)<\infty$.
\end{proposition}

{\sc Proof.} Let $L=V\begin{pmatrix}a&b\\c&d \end{pmatrix}$.  
First, consider an (infinite-dimensional) complement $X\subset V$ to $\ker a$ and consider
the subspace $L':=Xg$. We have $L'\cap W=0$.

 Consider the projection
$p(L')$ to $V$ along $W$. Consider a (finite-dimensional) complement $Z$ of $p(L')$ in $V$.
Set $L'':= L'\oplus Z$. We have $p(L'')=V$ and $L''\cap W=-0$. Thus
$L''\in \cM$.
\hfill $\square$

\sm

%%%%%%%%%%%%%%%%%%%%%%%%%%%%%%%%%%%%%%%%%%%%%%%%%%%%%%%%%%%%

{\bf\punct Atlas on $\bGr$.}
Consider finite subsets $\Omega$ and $\Xi$ in $\N$. Consider the following
subspaces $V[\Omega,\Xi]$ and $W[\Omega,\Xi]$ in $\ell\oplus\ell^\diamond$:
\begin{align}
V[\Omega,\Xi]:=\Bigl(\bigoplus\limits_{i\notin \Omega} \F_q e_i\Bigr)\oplus 
\Bigl(\bigoplus\limits_{j\in\Xi} \F_q f_j  \Bigr)
,
\label{eq:Vsk}
\\
W[\Omega,\Xi]:=\Bigl(\bigoplus\limits_{i\in \Omega} \F_q e_i\Bigr)\oplus 
\Bigl(\bigoplus\limits_{j\notin\Xi} \F_q f_j  \Bigr)
\label{eq:Wsk}
\end{align}
(sums are topological direct sums). We have
$$
\ell\oplus\ell^\diamond=V[\Omega,\Xi]\oplus W[\Omega,\Xi]
.
$$
Denote by $\cM[\Omega,\Xi]$ the set of all subspaces, which are graphs
of operators $V[\Omega,\Xi]\to W[\Omega,\Xi]$. 
Notice that 
$$
V=V[\varnothing,\varnothing],\qquad W=W[\varnothing,\varnothing],\qquad
\cM=\cM[\varnothing,\varnothing]
.$$

\begin{theorem}
\label{th:atlas}
$
\bGr=
\bigcup_{\Omega,\Xi\subset \Z} \cM[\Omega,\Xi]
.$ 
\end{theorem}

Consider a union $\N\sqcup\N'$ of two copies of $\N$, we assume that the first copy enumerates
vectors $e_j$ and the second copy vectors $f_i$.
Denote by $S(2\infty)$ the group of finite permutations $\N\sqcup\N'$. For each element
of $S(2\infty)$ we write 0-1 matrix 
(an element of $\GL(2\infty,\F_q)$) in the usual way.
We also add the following permutation to the group $S(2\infty)$:
$$
\dots\to e_2\to e_1\to f_1\to f_2\to \dots 
$$
It corresponds to the element $J\in \bGL(\ell\oplus\ell^\diamond)$.
Thus we get the extended group $\Z\ltimes S(2\infty)$.

\begin{lemma}
\label{l:symmetric}
An element $\sigma\in \Z\ltimes S(2\infty)\subset\bGL(\ell\oplus\ell^\diamond)$
 send a chart $\cM[\Omega,\Xi]$
to  the chart $\cM[\Omega_1,\Xi_1]$, where
$$
\bigl[(\N\setminus \Omega)\cup\Xi\bigr]\sigma = (\N\setminus \Omega_1)\cup\Xi_1
.$$
\end{lemma}

This is obvious.

\sm

{\sc Proof of Theorem \ref{th:atlas}.} The chart $\cM[\varnothing,\varnothing]$
is contained in $\bGr$. By Lemma \ref{l:symmetric} all charts $\cM[\Omega,\Xi]$
are contained in $\bGr$. 

\sm

Let us prove $\subset$. Let $L\in\cM[\Omega,\Xi]$. Let $M$ have codimension 1 in $L$.
We wish to show that $M$ is contained in some $\cM[\Omega',\Xi']$.
Consider the  projection $p(M)$ of $M$ to $V(\Omega,\Xi)$ along $W(\Omega,\Xi)$.
It is determined by an equation
$$
\sum_{j\notin \Omega} \alpha_j v_j + \sum_{i\in \Xi} \beta_i w_i=0,
\quad \text {where $\sum_{j\notin \Omega}  v_j e_j + \sum_{i\in \Xi}  w_i f_j\in L$}
$$
If some $\alpha_j\ne 0$, we can take $\Omega'=\Omega\cup\{j\}$, $\Xi'=\Xi$. 
If some $\beta_i\ne 0$, we can  take $\Omega'=\Omega$, $\Xi'=\Xi\setminus\{i\}$.

Next, let $L\in\cM[\Omega,\Xi]$ has codimension 1 in $K$. We wish to show
that $K$ is contained in some set $\cM[\Omega',\Xi']$.
We add a vector 
$$
h=\sum \beta_j e_j+\sum \gamma_i f_i
$$
to $L$. But this vector in not canonical, we can eliminate part of coordinates
and take
$$
h=\sum_{j\in \Omega} \gamma_j e_j+\sum_{i\notin\Xi} \delta_i f_i
.$$
If some $\gamma_j\ne 0$, we can exclude $\{j\}$ from $\Omega$. On the other hand,
if some $\beta_i\ne 0$, we can add $i$ to $\Xi$.

It remains to refer to Proposition
\ref{pr:LM}.
\hfill $\square$

%%%%%%%%%%%%%%%%%%%%%%%%%%%%%%%%%%%%%%%%%%%%%%%%%%%
%%%%%%%%%%%%%%%%%%%%%%%%%%%%%%%%%%%%%%%%%%%%%%%%%
%%%%%%%%%%%%%%%%%%%%%%%%%%%%%%%%%%%%%%%%%%%%%%%%%%%
%%%%%%%%%%%%%%%%%%%%%%%%%%%%%%%%%%%%%%%%%%%%%%%%%%

\section{The measure on Grassmannian}

\COUNTERS

{\bf \punct The measure $\mu$.} Note that each chart $\cM[\Omega,\Xi]$ is 
a direct product of countable family of  $\F_q$. We equip each chart by
the product measure.

\begin{theorem}
{\rm a)} There is a unique $\sigma$-finite measure $\mu$ on $\bGr$ coinciding
with the standard measure on each chart $\cM[\Omega,\Xi]$.
 
\sm

{\rm b)} This measure $\mu$ is $\bGL(\ell\oplus\ell^\diamond)$-invariant.

\sm

{\rm c)} The measure of $\bGr^0$ is
\begin{equation}
\mu(\bGr^0)=\prod_{j=1}^\infty(1-q^{-j})^{-1}.
\label{eq:total-measure}
\end{equation}

{\rm d)} Denote by $\cO_k$ the set of all $L\in\bGr^0$
such that $\dim L\cap W=k$. Then
\begin{equation}
\mu(\cO_k)= \frac{q^{-k^2}}{\prod_{j=1}^k (1-q^{-j})^2}
.
\label{eq:measure-orbit}
\end{equation}
 \end{theorem}  

The proof occupies the rest of this section.

\sm

%%%%%%%%%%%%%%%%%%%%%%%%%%%%%%%%%%%%%%%%%%%%%%%%%%%%%%%%%%%%%%%%%%%%% 

{\bf \punct The measure on $\ell^\diamond$.}
We equip $\F_q$ with uniform probability measure (the measure of each point is
$1/q$). We equip 
$$
\ell^\diamond=\F_q\times\F_q\times \F_q\times\dots
$$
with the product measure.
Note that this measure is the Haar measure on the compact Abelian group $\ell^\diamond$.
 We have natural maps $\pi_n:\ell^\diamond\to \F_q^n$
given by 
$$
(v_1,v_2,v_3\dots)\to(v_1,\dots,v_n)
.
$$
A  {\it cylindric subset} in $\ell^\diamond$ is a set which can be
obtained as   a preimage of a subset $S$ in $\F_q^n$ for some $n$. Finite unions and 
finite intersections of cylindric subsets are cylindric. 

\begin{proposition}
\label{pr:ellcirc-measure}
 The group $\bGL(\ell^\diamond)$ acts on $\ell^\diamond$
preserving the measure. Images and preimages of cylindric sets are cylindric.
\end{proposition}

This is obvious.

\sm

%%%%%%%%%%%%%%%%%%%%%%%%%%%%%%%%%%%%%%%%%%%%%%%%%%%%%%%%%%%%%%%%

{\bf \punct The measure on projective space.}
Consider  the projective space $\bbP\ell^\diamond=(\ell^\diamond\setminus 0)/\F_q^\times$.
Consider the map $\ell^\diamond\to\bbP\ell^\diamond$ defined a.s. Denote by 
$\kappa$ the image of the measure under this map. We have a quotient of an
invariant {\it finite} measure by an invariant partition 
($w\sim w'$ iff $w'=sw$, where $s\in\F_q^\times$). Therefore we get the following
trivial statement (it is used below for definition of the operator $\Delta$).

\begin{proposition}
\label{pr:projective-invariant}
The measure $\kappa$ on $\bbP\ell^\diamond$ is $\bGL(\ell^\diamond)$-invariant.
\end{proposition}

\sm

%%%%%%%%%%%%%%%%%%%%%%%%%%%%%%%%%%%%%%%%%%%%%%%%%%%%%%%%%

{\bf\punct The measure on $\cM$.} The set $\cM$ also is an infinite product 
of $\F_q$, we equip it with the product measure. Let us examine transformations
$$
\begin{pmatrix}a&b\\c&d \end{pmatrix}:
T\mapsto (a+Tc)^{-1}(b+Td)
$$
of $\cM$, see Lemma \ref{l:linear-fractional}.

\begin{lemma}
\label{l:parabolic-chart}
The group $\bfP$ acts on $\cM$ by measure preserving transformations.
They send cylindric sets to cylindric sets.
\end{lemma}

{\sc Proof.} Since $c=0$, we have transformations
$$
\begin{pmatrix}a&b\\0&d \end{pmatrix}:
T\mapsto a^{-1}b+ a^{-1}Td
$$
defined on the whole $\cM$. In fact we have the product of transformations of 3 types
corresponding to matrices $\begin{pmatrix}a&0\\0&\1 \end{pmatrix}$, 
$\begin{pmatrix}\1&0\\0&d \end{pmatrix}$,
$\begin{pmatrix}\1&b\\0&\1 \end{pmatrix}$.

The transformations $T\mapsto a^{-1}T$ act column-wise. Precisely,
denote by $t^{(1)}$, $ t^{(2)}$, \dots, the columns of $T$. Each column ranges in the space
$\ell^\diamond$, thus we can write
 $\cM=\ell^\diamond\times\ell^\diamond\times\dots$. Our transformation is
$$
(t^{(1)}, t^{(2)}, \dots)\to (a^{-1}t^{(1)}, a^{-1}t^{(2)}, \dots)
$$ 
and we refer to Proposition \ref{pr:ellcirc-measure}.

Similarly,  the transformations $T\mapsto Td$ act row-wise.

The transformations $T\mapsto T+b$ are simply translations on the compact Abelian group $\cM$.
\hfill $\square$.

\begin{lemma}
\label{l:cylindric}
Let $\begin{pmatrix}a&b\\c&d \end{pmatrix}\in\bGL^0(\ell\oplus\ell^\diamond)$.
Then the set $R$ of all $T\in\cM$ such that $\det(a+Tc)= 0$  is cylindric.
The map $T\mapsto(a+Tc)^{-1}(b+zd)$ is measure preserving outside $R$.
 \end{lemma} 

{\sc Proof.} First, let  $\begin{pmatrix}a&b\\c&d \end{pmatrix}\in \GL(2\infty,\F_q)$.
 We represent each block of $\begin{pmatrix}a&b\\c&d \end{pmatrix}$
as a block $(N+\infty)\times(N+\infty)$-matrix.
  For sufficiently large $N$ this matrix has the following structure
  $$
  \begin{pmatrix}
  a_{11}&0&b_{11}&0\\
  0&\1&0&0\\
  c_{11}&0&d_{11}&0\\
  0&0&0&\1
  \end{pmatrix}
  .
  $$ 
Then our transformation has the form
{\small
\begin{multline*}
\begin{pmatrix} T_{11}& T_{12}\\T_{21}&T_{22}  \end{pmatrix}
\mapsto
\\
\left[\begin{pmatrix} a_{11}&0\\0&\1\end{pmatrix}+
\begin{pmatrix} T_{11}& T_{12}\\T_{21}&T_{22}  \end{pmatrix}
 \begin{pmatrix} c_{11}&0\\0&0\end{pmatrix}\right]^{-1}
 \left[
  \begin{pmatrix} b_{11}&0\\0&0\end{pmatrix}+
 \begin{pmatrix} T_{11}& T_{12}\\T_{21}&T_{22}  \end{pmatrix}
  \begin{pmatrix} d_{11}&0\\0&\1\end{pmatrix}
 \right]
=\\=
\begin{pmatrix}a_{11}+T_{11}c_{11}&0\\ 
T_{21}c_{11}&\1
  \end{pmatrix}^{-1} 
\begin{pmatrix}b_{11}+T_{11}d_{11}&T_{12}\\ 
T_{21}d_{11}&T_{22}
  \end{pmatrix}=\\=
 \begin{pmatrix}(a_{11}+T_{11}c_{11})^{-1}&0\\ 
-T_{21}c_{11}(a_{11}+T_{11}c_{11})^{-1}&\1
  \end{pmatrix}
  \begin{pmatrix}b_{11}+T_{11}d_{11}&T_{12}\\ 
T_{21}d_{11}&T_{22}
  \end{pmatrix}
  =\\%=
  \begin{pmatrix}
  (a_{11}+T_{11}c_{11})^{-1}(b_{11}+T_{11}d_{11})&(a_{11}+T_{11}c_{11})^{-1}T_{12}
  \\
 T_{21}\!\!\left[-c_{11}(a_{11}+T_{11}c_{11})^{-1}(b_{11}+T_{11}d_{11})+d_{11}\right]
 &-T_{21}c_{11}(a_{11}+T_{11}c_{11})^{-1}T_{12}+T_{22}
  \end{pmatrix}
\end{multline*}
}
We observe that the set $R$ is determined by the equation
$\det(a+T_{11}d_{11})=0$. Evidently, this set is cylindric. Fix $T_{11}$ outside
'surface' $\det(a+T_{11}d_{11})=0$. The set $T_{11}=const$ is a product of 3
spaces $\F_q^{N\times\infty}$, $\F_q^{\infty\times N}$, $\F_q^{\infty\times\infty}$
corresponding blocks $T_{12}$, $T_{21}$, $T_{22}$ of $T$.

The transformations $T_{12}\mapsto(a_{11}+T_{11}c_{11})^{-1}T_{12}$ are column-wise
(columns have length $N$).

The transformations 
$T_{21}\mapsto  T_{21}\left[-c_{11}(a_{11}+T_{11}c_{11})^{-1}(b_{11}+T_{11}d_{11})+d_{11}\right]$ 
are row-wise. Rows have length $N$. The matrix in square brackets 
is invertible. Indeed,
\begin{multline*}
0\ne
\det \begin{pmatrix}a_{11}&b_{11}\\c_{11}&d_{11} \end{pmatrix}
=
\det \begin{pmatrix}a_{11}+T_{11}c_{11}&b_{11}\\c_{11}+T_{11}d_{11}&d_{11} \end{pmatrix}
=\\=
\det(a+zc)\det\left[ -c_{11}(a_{11}+T_{11}c_{11})^{-1}(b_{11}+T_{11}d_{11})+d_{11}\right]
\end{multline*}

It remains to examine transformations 
$$T_{22}\to -T_{21}c_{11}(a_{11}+T_{11}c_{11})^{-1}T_{12}+T_{22}.$$
First, let $k,l>N$. The new matrix element $t_{kl}$ of $T$ depends only 
on elements $t_{\alpha,\beta}$ with $\alpha\le k$, $\beta\le l$. Therefore the preimage of any
 cylindric set
is cylindric. We can apply the same argument to the inverse transformation
and get that the image of any cylindric set is cylindric.

  Next, for fixed $T_{12}$, $T_{21}$ the transformation of $T_{22}$ is a translation. Therefore it preserves the measure.
  
  Now let $\begin{pmatrix}a&b\\c&d \end{pmatrix}\in\bGL^0(\ell\oplus\ell^\diamond)$
  be arbitrary. We refer to Proposition \ref{pr:cor}, recall that
   each element of $\bGL^0$ can be represented
  as
  $$
  g=h q  r, \qquad\text{where $h\in\bGL(\ell)$, $p\in\GL(2\infty,\F_q)$, $r\in\bfP$}.  
  $$
 By Lemma \ref{l:parabolic-chart} $h$ and $r$ act on $\cM$ by measure-preserving
 transformations defined everywhere, hence the lemma reduces to the case
 examined above.
   \hfill $\square$
  
\sm

Each chart $\cM[\Omega,\Xi]$ is the  space of square matrices, we equip a chart
$\cM[\Omega,\Xi]$ by
 the product measure $\mu_{\Omega,\Xi}$.

\begin{lemma}
\label{l:coincide-intersection}
The measures $\mu_{\Omega,\Xi}$ and $\mu_{\Omega',\Xi'}$
coincide on $\cM[\Omega,\Xi]\cap\cM[\Omega',\Xi']$.
\end{lemma}

{\sc Proof.} Without loss of generality, we can assume 
$\cM[\Omega',\Xi']=\cM[\varnothing, \varnothing]=\cM$. The intersection
$\cM\cap \cM[\Omega,\Xi]$ is non-empty if and only if $\#\Omega=\#\Xi$. Consider an element
$\sigma\in S(2\infty)\subset \bGL^0$ sending $\cM$
 to $\cM[\Omega,\Xi]$,
see Lemma \ref{l:symmetric}.  The map $\sigma$ send the measure $\mu_{\varnothing, \varnothing}$
to the  measure $\mu_{\Omega,\Xi}$. Represent $\sigma$, $\sigma^{-1}\in \bGL^0$ in the block
form 
$$\sigma=\begin{pmatrix}\alpha&\beta\\ \gamma&\delta\end{pmatrix},
\qquad
\sigma^{-1}=\begin{pmatrix}\alpha^t&\gamma^t\\\beta^t& \delta^t\end{pmatrix}
,
$$
where $^t$ denotes transposition.
Let $L\in \cM$  be a graph of $S:\ell\to\ell^\diamond$.
Applying $\sigma$ to $L$ we get a subspace, whose coordinate in $\cM$ 
is $T= (\alpha+S\gamma)^{-1}(\beta+S\delta)$, see Lemma \ref{l:linear-fractional}.d. In fact, we get a bijection
\begin{equation}
\cM\cap \sigma^{-1}\cM [\Omega,\Xi]
\,\,
\to\,\, \cM\cap \cM [\Omega,\Xi]
\end{equation}
In coordinates, the domain of definiteness of the map consists
of $S$ such that $\det(\alpha+S\gamma)\ne 0$. The image consists of $T$ such that
$\det(\alpha^t+\beta^t T)\ne 0$, since the inverse map is given by 
$S=(\alpha^t +T\beta^t)^{-1}(\gamma^t+S\delta^t)$. 
By Lemma \ref{l:cylindric}, a linear fractional map preserves measure $\mu$
on the domain of definiteness. \hfill $\square$

\begin{corollary}
The measure $\mu$ is well-defined.
\end{corollary}

Next,  an element $g\in\bGL^0(\ell\oplus\ell^\diamond)$
determines a partially defined map $\cM[\varnothing, \varnothing]\to \cM[\Omega,\Xi]$
for any $\cM[\Omega,\Xi]$.

\begin{lemma}
This map is measure preserving.
\end{lemma}

{\sc Proof.} Indeed, according Lemma \ref{l:linear-fractional}, 
in coordinates  this
map  is given by a  linear-fractional expression (as in the proof  of
Lemma \ref{l:coincide-intersection}). \hfill $\square$

\sm

We also apply the same argument for arbitrary pair of charts of $\bGr^\alpha$ and get the 
the $\bGL^0$-invariance of the measure $\mu$ on $\bGr^\alpha$. The invariance of $\mu$
 with respect to the
operator $J$ is evident.

%\begin{corollary}
%The measure $\mu$ is $\bGL(\ell\oplus\ell^\diamond)$-invariant.
%\end{corollary}

%%%%%%%%%%%%%%%%%%%%%%%%%%%%%%%%%%%%%%%%%538

{\bf\punct Uniqueness of $\mu$.} We have a natural Borel structure
on each chart and therefore a natural Borel structure on $\bGr$.

\begin{proposition}
\label{pr:measure-unique}
{\rm a)} A $\bGL^0(\ell\oplus \ell^\diamond)$-invariant finite or $\sigma$-finite
measure on $\bGr^j$  defined on all Borel sets coincides with $\mu$ up to a constant factor.

\sm

{\rm b)} A  $\bGL(\ell\oplus\ell^\diamond)$-invariant $\sigma$-finite
 measure  on $\bGr$  defined on all Borel sets
coincides with $\mu$ up to a constant factor.
\end{proposition}

{\sc Proof.} 
The restriction of our measure
 to $\cM\simeq \F_q^{\infty\times\infty}$ must be invariant with respect
to all translations $T\mapsto T+A$, where $A$ is an operator $\ell\to \ell^\diamond$.
 Therefore we get the Haar measure on the compact group $\F_q^{\infty\times\infty}$
 (i.e., our measure $\mu$ on $\cM$).
 For any chart $\cM[\Omega,\Xi]$ we take a map sending $\cM\to \cM[\Omega,\Xi]$. This
 fixes measures on all charts.  \hfill $\square$

%%%%%%%%%%%%%%%%%%%%%%%%%%%%%%%%%%%%%%%%%%%%%%%%%

\sm

{\bf \punct Finite Grassmannians.} Denote by $\Gr_m^k$ the set of all
$k$-dimensional subspaces in the linear space $\F_q^m$.
Let $\Gr_m$ be set of all subspaces in $\F_q^m$,
$$
\Gr_m=\bigsqcup_{j=0}^m \Gr_m^k
.$$
We need some simple facts (see, e.g. \cite{Cecc}) about 
such Grassmannians.

\sm

1) Number of elements in $\GL(m,\F_q)$ is
$$
\gamma_m:=
\prod_{j=1}^m (q^m-q^j)=q^{m^2} \prod_{j=1}^m (1-q^{-j}).
$$

2) Denote by $P=P_n\subset \GL(2n,\F_q)$ the stabilizer of the subspace 
$0\oplus\F_q^n$ in $\F_q^{2n}$.
 The number of elements of $\Gr_{2n}^n$
is 
\begin{equation}
\frac{\#\GL(2n,\F_q)}{\# P}=
\frac{\gamma_{2n}}{\gamma_n^2\, q^{n^2}}=q^{n^2}\cdot
\frac{\prod_{j=1}^{n} (1-q^{-n-j} )}{\prod_{j=1}^{n}(1-q^{-j})}
\label{eq:gr-2n}
.\end{equation}
Indeed, $\Gr_{2n}^n$ is a $\GL(2n,\F_q)$-homogeneous space,
the subgroup $P$ consists of $(n+n)\times(n+n)$-matrices
$\begin{pmatrix}a&b\\0&d\end{pmatrix}$, where $a$, $d\in\GL(n,\F_q)$
and $b$ is arbitrary.

\sm

3) Consider the set $\cO_k(n)$ of
 all subspaces $L\in \Gr_{2n}^n$ that have $k$-dimensional
intersection with $W:=0\oplus \F_q^n$. The number of elements of this set is
\begin{equation}
\#\cO_k(n)=
\frac{\gamma_n^2\, q^{n^2}}
{\gamma_n^2 \gamma_{n-k}^2\, q^{4k(n-k)}}=q^{n^2}\cdot
\frac{q^{-k^2}\prod_{j=n-k+1}^n(1-q^{-j})^2} {\prod_{j=1}^k(1-q^{-j})^2}
\label{eq:number-intersection}
\end{equation}
Indeed, the  stabilizer $Q$ of a pair  $(L,M)$ of subspaces stabilizes also $L\cap M$, $L+M$
and the decomposition of $(L+M)/(L\cap M)$ into a direct sum
$$(L+M)/(L\cap M)= L/(L\cap M)\oplus M/(L\cap M)$$
Therefore $Q$ consists of invertible
 matrices of size $k+(n-k)+(n-k)+k$ having the following structure
$$
\begin{pmatrix}*&*&*&*\\
0&*&0&*\\
0&0&*&*\\
0&0&0&*
 \end{pmatrix}
$$
The desired number is $\frac{\#P}{\#Q}$. \hfill $\square$

\sm

4) We regard $\F_q^{2n}$ as a direct sum $V_n\oplus W_n:=\F_q^n\oplus\F_q^n$, set
$$
V:=\bigoplus_{j=1}^N \F_q e_j,\qquad W:=\bigoplus_{j=1}^N \F_q f_j
.
$$
Next, for subsets $\Omega$, $\Xi\in\{1,2,\dots,n\}$ such that $\#\Omega=\#\Xi$,
we define subspaces $V_n[\Omega,\Xi]$, $W_n[\Omega,\Xi]$ as above, 
(\ref{eq:Vsk})--(\ref{eq:Wsk}).
As above, we define a chart $\cM_n[\Omega,\Xi]\subset \Gr_{2n}^n$ as the set of graphs 
of operators $V_n[\Omega,\Xi]\to W_n[\Omega,\Xi]$. 

We define a uniform measure $\mu_n$ on $\Gr_{2n}^n$ by the assumption:
for any set $S\subset\Gr_{2n}^n$,
$$
\mu_n(S)=\frac{\#S}{q^{n^2}}
.
$$
In particular, the measure of any chart $\cM_n[\Omega,\Xi]$ is 1.

%%%%%%%%%%%%%%%%%%%%%%%%%%%%%%%%%%%%%%%%%%%%%%%%%%%%%%%%%

\sm

{\bf\punct The map $\bGr^0\to\Gr_{2n}$.%
\label{ss:pi}} Consider the the following subspaces
$X_n$ and $Y_n$ in $\ell\oplus\ell^\diamond$:
$$
X_n:=\bigoplus_{j>n} \F_q f_j\qquad
Y_n:=\bigoplus_{i\le n}\F_q e_i \oplus \bigoplus_{j\in \N} \F_q f_j
.
$$
We have an obvious isomorphism 
$$
Y_n/X_n\simeq \F_q^{2n}.
$$
For $L\in\bGr^0$ we define $\pi_n(L)\in \Gr_{2n}$
as the image of $L\cap Y_n$ under the map $Y_n\to Y_n/X_n$.

\sm

{\sc Remark.} A dimension of $\pi_n(L)$ can be arbitrary between 0 and $2n$.
%\hfill$\square$

\begin{lemma}
{\rm a)}
If $\Omega$, $\Xi\subset\{1,2,\dots,n\}$,  and $L\in\cM[\Omega,\Xi]$,
 then $\dim \pi_n(L)=n$.

\sm

{\rm b)}
Moreover, under the same condition
$$
\pi_n(\cM[\Omega,\Xi])= \cM_n[\Omega,\Xi]
$$
and the image of the measure $\mu$ on $\cM[\Omega,\Xi]$
is the measure $\mu_n$ on $\cM_n[\Omega,\Xi]$.

\sm

{\rm c)}  Any cylindric subset in $\cM[\Omega,\Xi]$
is a $\pi_n$-preimage of a subset in $\cM_n[\Omega,\Xi]$
for sufficiently large $n$.
\end{lemma}

This is straightforward.

\sm

Denote by $\cU_n$ the following subset in $\bGr^0$:
$$
\cU_n=\bigcup_{\begin{array}{c}
\Omega, \Xi\subset\{1,2,\dots,n\}\\ \#\Omega=\#\Xi
\end{array}}
\cM[\Omega,\Xi]
$$
Let $\phi$ be a function on $\Gr_{2n}^n$.
Consider the function  $\phi^*$ on $\bGr^0$ given by
\begin{equation}
\phi^*(L):=
\begin{cases}
\phi(\pi_n(L)),\qquad& \text{if $L\in\cU_n$;}
\\
0 &\text{otherwise}
\end{cases}
\label{eq:phi*}
\end{equation}

{\sc Remark.} The preimage of $\Gr_{2n}^n$ under the map 
$\pi_n$ is larger than $\cU_n$.
%\hfill $\square$

\sm

Denote by $\cF_n$ the set of functions on $\bGr^0$
that can be obtained in this way. Obviously,
$\cF_{n+1}\supset \cF_n$.

\begin{lemma}
\label{l:union}
{\rm a)} The map $\pi_n:\cU_n\to\Gr_{2n}^n$ commutes with action of $\GL(2n,\F_q)$.

\sm

{\rm b)} $\bigcup_n \cU_n=\bGr^0$. 

\sm

{\rm c)} $\bigcup_n \cF_n$ is dense in $L^2(\bGr^0)$.
\end{lemma}

This is obvious.

\begin{corollary}
The total measure of $\bGr^0$ is finite and equals 
{\rm(\ref{eq:total-measure})}.
\end{corollary}

{\sc Proof.}
$$
\mu(\bGr^0)=\lim_{n\to\infty} \cU_n=
\lim_{n\to\infty} \mu_n(\Gr_{2n}^n),
$$
and we apply (\ref{eq:gr-2n}). \hfill $\square$

\begin{corollary}
The measure of $\cO_k$ is given by  {\rm (\ref{eq:measure-orbit})}.
\end{corollary}

{\sc Proof.} Let $L\in \bGr^0$. We have $W\subset Y_n$.
On the other hand for large $n$, we have $L\cap X_n=0$ (since $L\cap W$
is finite-dimensional). Therefore, starting sufficiently large $n$ we have
$$
\dim (\pi_n(L)\cap W_n)=\dim L\cap W,
$$
Therefore
$$
\mu(\cO_k)=\lim_{n\to\infty} \cO_k(n)
$$
and we apply (\ref{eq:number-intersection}).

%%%%%%%%%%%%%%%%%%%%%%%%%%%%%%%%%%%%%%%%%%%%%%%%%%%
%%%%%%%%%%%%%%%%%%%%%%%%%%%%%%%%%%%%%%%%%%%%%%%%%%%%%
%%%%%%%%%%%%%%%%%%%%%%%%%%%%%%%%%%%%%%%%%%%%%%%%%%%%%
%%%%%%%%%%%%%%%%%%%%%%%%%%%%%%%%%%%%%%%%%%%%%%%%%%%%%%

\section{Decomposition of $L^2(\bGr^0)$ and
\\ Al-Salam--Carlitz polynomials}

\COUNTERS

The group $\bGL^0(\ell\oplus\ell^\diamond)$ acts in $L^2(\bGr^0)$ by changes of variables
$$
\rho(g)F(L)=F(Lg).
$$
Here we decompose this representation.

\sm

{\bf \punct Operator of average.}
Let $L\in\bGr$. Denote by

\sm

--- $\Sigma^\downarrow_L$ the set of all subspaces $K\subset L$ such that
$\dim (L/K)=1$; 

\sm

--- $\Sigma^\uparrow_L$ the set of all spaces $N\in\bGr$ such that $N\supset L$
and $\dim N/L=1$;

\sm

--- $\Sigma^\updownarrow_L$ the set of subspaces $M\in\bGr$ such that 
$$
\dim L/(L\cap M)=1,\qquad \dim M/(L\cap M)=1.
$$

\begin{lemma}
\label{l:nu}
 Let $L\in\bGr$, $\cP_L\subset \bGL(\ell\oplus\ell^\diamond)$ be the stabilizer of $L$ in $\bGL(\ell\oplus\ell^\diamond)$.
 
\sm 
 
{\rm a)} There is a unique  $\cP(L)$-invariant probability measure $\nu^\updownarrow_L$ on 
the space $\Sigma^\updownarrow_L$.

\sm

{\rm b)} There is a unique  $\cP(L)$-invariant probability measure $\nu^\downarrow_L$ on 
the space $\Sigma^\downarrow_L$. 

\sm

{\rm c)} There is a unique  $\cP(L)$-invariant probability measure $\nu^\uparrow_L$ on 
the space $\Sigma^\uparrow_L$. 
\end{lemma}

{\sc Proof.} b) Since the group $\bGL$ acts on $\bGr$ transitively,
without loss of generality we set $L=V$. The group
$\cP$ consists of matrices $\begin{pmatrix}a&0\\c&d\end{pmatrix} $.
Subspaces of codimension 1 in $V$ are determined by equations $\sum \alpha_j v_j=0$,
where $\sum v_j e_j\in V$ and
$\alpha=(\alpha_1,\alpha_2,\dots)$ is in $\ell^\diamond$. Therefore the set of subspaces $K$
is the projective space $\bbP\F_q^\infty\simeq(\ell^\diamond\setminus 0)/\F_q^\times$
and equipped with a canonical measure $\kappa$, see Proposition \ref{pr:projective-invariant}.
 Matrices 
$\begin{pmatrix}1&0\\c&d\end{pmatrix} $ fix all elements of $V$ and therefore act trivially
on the projective space $\bbP\ell^\diamond$. The group of matrices
$\begin{pmatrix}a&0\\0&\1\end{pmatrix} $ acts by projective transformations and therefore
preserve the measure $\kappa$. 

\sm

c) Again, without loss of generality we set $L=V$. The quotient space 
$(\ell\oplus\ell^\diamond)/V$ is isomorphic to $\ell^\diamond$. Therefore
overspaces $M\supset L$ are enumerated by points of $\bbP\ell^\diamond$.

\sm

a) The existence of the measure follows from b) and c). The space $\Sigma^\updownarrow_L$
is $\cP$-homogeneous, therefore the invariant measure is unique.
\hfill $\square$ 

\sm

We define the operator $\Delta$ in $L^2(\bGr^0)$
by
$$
\Delta f(L)=\int_{M\in\Sigma^\updownarrow_L} f(M)\,d\nu^\updownarrow_L(M)
.
$$

\begin{proposition}
\label{pr:Delta}
{\rm a)} The operator $\Delta$ is a self-adjoint bounded operator 
in $L^2(\bGr^0)$.

\sm

{\rm b)} The operator $\Delta$ is $\bGL^0$-intertwining.
\end{proposition}

{\sc Proof.} The statement b) follows from the $\cP_L$-invariance
of the measure $\nu_L$.

Consider the measure $\nu^\updownarrow$ on $\bGr^0\times \bGr^0$ defined by
$$
\nu^\updownarrow(S)=\int_{\bGr^0} \nu^\updownarrow_L (S\cap (\{L\}\times \bGr^0)\,d\mu(L),
\qquad S\subset \bGr^0\times \bGr^0
,$$
here we regard the  measure $\{\nu^\updownarrow_L\}$
 as a measure on $\bGr^0$. We identify
 any fiber $\{L\}\times \bGr^0\subset \bGr^0\times \bGr^0$ with $\bGr^0$.
This measure satisfies
\begin{equation}
\la f, \Delta g\ra_{L^2(\bGr^0)}=
\int_{\bGr^0\times \bGr^0} f(L)\,\ov g(M) \,d\nu^\updownarrow(L,M)
\label{eq:markov}
\end{equation}

We need the following lemma.

\begin{lemma}
\label{l:nu-symmetry}
The measure $\nu^\updownarrow$ is symmetric with respect to transposition of factors
in $\bGr^0\times \bGr^0$.
\end{lemma}

{\sc Proof of Lemma.}
The measure $\nu^\updownarrow$ is invariant with respect to the diagonal action
of $\bGL^0(\ell\oplus\ell^\diamond)$ on  $\bGr^0\times\bGr^0$.
Therefore its projection to the second factor is  $\bGr^0$-invariant,
 by Proposition \ref{pr:measure-unique} it coincides with
 $\mu$. Consider conditional measures on the fibers $\bGr^0\times \{K\}$.
 They must be  $\cP_L$-invariant a.s., hence they coincide with $\nu^\updownarrow_L$.
  This implies the desired statement.   
\hfill $\square$

\sm

{\sc End of proof of Proposition \ref{pr:Delta}.}
Let us look to (\ref{eq:markov}) as to an abstract expression. Let  $X$ be a space with
a probability measure
$\mu$, let $S$ be a measure on the product $X\times X$ and projections of $S$ to both factors coincide with 
$\mu$. Then the formula (\ref{eq:markov}) determines a so-called 
  {\it Markov operator}. Such operators  automatically satisfy
$\|\Delta\|\le 1$, see, e.g., \cite{Ner-book}, Theorem VIII.4.2.
Since the measure $S$ is symmetric with respect to the transposition of factors, 
$\Delta$ is self-adjoint.
\hfill $\square$
 
%%%%%%%%%%%%%%%%%%%%%%%%%%%%%%%%%%%%%%%%%%%%%%%   

\sm

{\bf \punct $\bfP$-invariant functions.} Recall that $\bfP$ is the stabilizer
of the subspace $W=0\oplus \ell^\diamond$ in $\bGL^0(\ell\oplus\ell^\diamond)$. Orbits $\cO_k$
of $\bfP$ on the Grassmannian $\bGr^0$ are the following sets:
$$
\cO_k=\{
L:\quad \dim( L\cap W) =k\}.
$$

Denote by $\cH$ the space of $\bfP$-fixed functions on $\bGr^0$,
it consists of functions constant on orbits $\cO_k$.

\begin{proposition}
\label{pr:P-fixed}
The $\bGL^0(\ell\oplus\ell^\diamond)$-cyclic span of $\cH$
is the whole $L^2(\bGr^0)$.
\end{proposition}

A proof is contained in the next subsection.

\begin{corollary}
Any $\bGL^0$-invariant subspace $V$ of $L^2(\bGr^0)$ is a $\bGL^0$-cyclic span
of $\bfP$-fixed vectors contained in $V$
\end{corollary}

This statement implies Theorem 1.5.a.

\sm

{\sc Proof of the corollary.} Denote by $W$ the cyclic span of $\bfP$-fixed vectors contained
in $V$. Let $Y$ be the orthocomplement of $W$ in $V$, let $Z$ be the orthocomplement of
$Y$ in the whole $L^2$. Then all $\bfP$-fixed vectors are contained in $Z$ and therefore their
cyclic span also is contained in $Z$. Therefore $Y=0$.
\hfill $\square$

\sm

{\sc Remark.}  The representation of $\bGL^0$ in $L^2(\bGr^0)$ admits a single cyclic 
$\bfP$-invariant vector. Equivalently, $L^2$ is a direct sum of pairwise
nonequivalent subrepresentations.  We omit a  proof of this statement (this requires additional calculations).

\sm

{\bf\punct Proof of Proposition \ref{pr:P-fixed}.}
Consider the group $\GL(2n,\F_q)$ and the homogeneous
 space $\Gr_{2n}^n=\GL(2n,\F_q)/P$, where $P=P_n$ is
 the group of $(n+n)\times(n+n)$-matrices 
 $\begin{pmatrix}a&b\\0&d \end{pmatrix}$. 
 The pair $\GL(2n,\F_q)\supset P$ is spherical
 (see, e.g. \cite{Cecc},  8.6.5). The representation
 of $\GL(2n,\F_q)$ on the space of functions on $\Gr_{2n}^n$ was 
 decomposed in \cite{ZZ}, \cite{Sta}. It is a sum of $n+1$
 irreducible subrepresentations, each subrepresentation
 has a unique $P$-fixed vector (by the Frobenius reciprocity).
 Such vectors can be regarded as functions of $k=\dim L\cap W$. 
  $P$-fixed vectors are given by
 $q$-Hahn polynomials%
 \footnote{On $q$-Hahn polynomials, see, e.g., \cite{Koe}.}
 $Q_0$, $Q_1$, \dots, $Q_n$:
 $$
 Q_j(q^{-k};q^{-n-1},q^{n-1};n;q):=
 {}_3\phi_2\left(
 \begin{matrix}q^{-j},q^{j-2n-1},q^{-k}\\q^{-n},q^{-n} \end{matrix};q;q
 \right)
 .$$
$q$-Hahn polynomials are regarded as functions of the discrete variable $k=0$, 1,\dots, $n$.
 
 We use the standard notation for basic hypergeometric functions
 \begin{multline}
  {}_r\phi_s\left(\begin{matrix}a_1,\dots,a_r\\b_1,\dots,b_s\end{matrix};q,z  \right)
  :=\\:=
  \sum_{k=0}^\infty \frac{(a_1;q)_k \dots (a_r;q)_k}{(b_1;q)_k(b_s;q)_k}
  (-1)^{(1+s-r)k} q^{(1+s-r)k(k-1)/2} \frac{z^k}{(q;q)_k}
,\end{multline}
and
$$
(a;q)_k:=(1-a)(1-aq)\dots(1-aq^{k-1}).
$$
Consider the difference operator
\begin{equation}
\cL y(k)=B(k) y(k+1)-(B(k)+D(k))y(k)+D(k)y(k-1),
\label{eq:L-hahn}
\end{equation}
where
\begin{align}
B(k)&=(1-q^{k-n})^2;
\label{eq:B(k)}
\\
D(k)&=(1-q^{-2n-2})(1-q^k)^2
\label{eq:D(k)}
\end{align}
The $q$-Hahn polynomials are eigenfunctions of $\cL$:
$$
\cL Q_j(k)=(1-q^{-j})(1-q^{j-2n-1})\, Q_j(k)
.$$

Consider the function $I_n(L)$ on $\Gr_{2n}^n$,  which equals
$1$ on $\cM_n[\varnothing,\varnothing]$ and 0 otherwise

\begin{lemma}
\label{l:cyclic}
The $\GL(2n,\F_q)$-cyclic span of $I_n$ is the whole space of functions
on $\Gr_{2n}^n$. 
\end{lemma}

{\sc Proof.} Assume that there is an irreducible subrepresentation $Y$
orthogonal to $I_n$. It contains a $P$-fixed function $Q_j(k)$, which also must
be orthogonal to $I_n$. This means that  $Q_j(0)=0$.

By (\ref{eq:B(k)}), $B_j(k)$ is non-zero for all $k=0$, 1,\dots, $n-1$. By
(\ref{eq:L-hahn}) we can express
\begin{align*}
y_j(k+1)&=s_k\cdot y(k)+t_k\cdot y(k-1);
\\
y(1)&= s_0 y(0)
,\end{align*}
where $s_k$, $t_k$ are some coefficients. Therefore $Q_j(0)=0$
implies $Q_j(k)=0$ for all $k$.
\hfill $\square$

\sm

%%%%%%%%%%%%%%%%%%%%%%%%%%%%%%%%%%%%%%%%%%%%%%%%%%%%

{\sc  Proof of Proposition \ref{pr:P-fixed}.} Now consider the characteristic
function $I_\infty$ of the chart $\cM[\varnothing,\varnothing]\subset \Gr_{2n}^n$. 
This function is contained in the subspace $\cF_n$ (see Subsection \ref{ss:pi})
and has the form $I_n^*(L)$, see (\ref{eq:phi*}). By Lemma \ref{l:cyclic}
the cyclic span of $I_n$ is the whole space of functions on Grassmannian.
By Lemma \ref{l:union}.a, the $\GL(2n,\F_q)$-cyclic span of $I_\infty$ 
is $\cF_n$.  Therefore, $\GL(2\infty,\F_q)$-cyclic span of $I_\infty$
is $\bigcup_n \cF_n$. The latter subspace is dense in $L^2(\bGr^0)$,
see Lemma \ref{l:union}.d.
\hfill $\square$

\sm

%%%%%%%%%%%%%%%%%%%%%%%%%%%%%%%%%%%%%%%%%%%

{\bf \punct Intersections of $\Sigma_L$ with $P$-orbits.}
Let $L\in\cO_k$. The set $\Sigma_L$ have intersections with $\cO_{k-1}$, $\cO_k$,
$\cO_{k+1}$. Let us find $\nu^\updownarrow_L$-measures of these intersections.

\begin{lemma}
\label{l:jumps}
\begin{align}
 \nu^\updownarrow_L(\Sigma_L\cap \cO_{k-1})&= (1-q^{-k})^2;
 \label{eq:k-1}
 \\
  \nu^\updownarrow_L(\Sigma_L\cap \cO_{k})&=q^{-k}-q^{-2k}-q^{-2k-1};
  \label{eq:k}
  \\
   \nu^\updownarrow_L(\Sigma_L\cap \cO_{k+1})&= q^{-2k-1}.
   \label{eq:k+1}
\end{align}
\end{lemma}

{\sc Proof.} We will choose a subspace $K\subset L$
of codimension 1. Next, we choose an overspace $M\supset K$
such that $\dim M/K=1$.

{\sc Step 1.}
Without loss of generality we can assume
$$
L=\bigoplus_{j>k}\F_q e_j\oplus \bigoplus_{i\le k} \F_q f_i 
$$
and 
\begin{equation}
\label{eq:LW}
L\cap W=\bigoplus_{i\le k} \F_q f_i.
\end{equation}
A subspace $M$ is determined by a linear equation
$$
\sum_{j>k} \alpha_j v_j+ \sum_{i\le k} \beta_i w_i=0,
\quad\text{where $\sum_{j>k}  v_j e_j+ \sum_{i\le k}  w_i f_j\in V$}
.$$

{\sc Case 1.} All $\beta_i=0$. Then $M\cap W=L\cap W$. 
The probability of this event is $q^{-k}$.

\sm

{\sc Case 2.} There exists $\beta_i\ne 0$. Then $W\cap M$ is the subspace
in (\ref{eq:LW}) determined by the equation $\sum_{i\le k} \beta_i w_i=0$.
Thus $\dim M\cap L=k-1$. The probability of this event is $1-q^{-k}$.

\sm

{\sc Step 2.}
Next, we choose $K\supset M$.

In Case-1 we by rotation set the subspace $K$ to the position
$$
K=\bigoplus_{j>k+1}\F_q e_j\oplus \bigoplus_{i\le k} \F_q f_i 
.
$$
We add a vector
$$
h=\sum_{j\le k+1} x_j e_j+ \sum_{j>k} y_j f_j
$$
to $K$. Again, there are two cases.

\sm

{\sc Case 1.1.}  All $x_j=0$. Then $M\cap W=\F_q h\oplus \bigoplus_{i\le k} \F_q f_i$.
Therefore $M\in \cO_{k+1}$. Conditional probability of this event is $q^{-k-1}$. 
 
\sm 
 
 {\sc Case 1.2.} There exists $x_j\ne 0$. Then $M\cap W=K\cap W$
 and $M$ falls to $\cO_k$.
 The conditional probability of this event is $1-q^{-k-1}$.

\sm

In Case-2 we put the subspace $K$ to the position
$$
K=\bigoplus_{j>k}\F_q e_j\oplus \bigoplus_{i\le k-1} \F_q f_i 
,
$$
add a vector
$$
h=\sum_{j\le k} x_j e_j+ \sum_{j>k-1} y_j f_j
$$
to $K$. Repeating the same consideration as in
Cases 1.1-1.2 we get

{\sc Case 2.1.} With conditional probability  $q^{-k}$ we fall to $\cO_{k}$.

{\sc Case 2.2.} With conditional probability  $(1-q^{-k})$ we fall to $\cO_{k-1}$.

Compiling 4 sub-cases we get (\ref{eq:k-1})--(\ref{eq:k+1}).
\hfill $\square$

\sm

%%%%%%%%%%%%%%%%%%%%%%%%%%%%%%%%%%%%%%%%%%%%%%%%%%%%%%%%%%%%%%%

{\bf\punct Difference operator.}
Now let us restrict the operator $\Delta$ to the subspace $\cH$
of $\bfP$-invariant functions. The space $\cH$ can be regarded as 
the space of complex sequences
$$
x=(x(0), x(1), x(2),\dots)
$$
with inner product
$$
\la x,y\ra=\sum_{j=0}^\infty \frac{q^{-k^2}}{\prod_{j=1}^k(1-q^{-j})^2}\cdot x(k)\ov y(k)
.
$$

\begin{proposition}
The restriction of $\Delta$ to $\cH$ is the following difference operator
$$
\cL y(k)=(1-q^{-k})^2 y(k-1)+(2q^{-k}-q^{-2k}-q^{-2k-1}) y(k)+q^{-2k+1} y(k+1)
$$
\end{proposition}

{\sc Proof.} The statement follow from Lemma \ref{l:jumps}.
\hfill $\square$

\sm

This difference Sturm--Liouville problem is known, the solutions are 
Al-Salam--Carlitz II polynomials
$$
V_j^{(1)}(q^k,q^{-1})=(-1)^j q^{j(j-1)/2}
{}_2\phi_0\left( \begin{matrix} q^j, q^k\\ -\!- \end{matrix};q^{-1};q^{-j}\right)
$$
the corresponding eigenvalues are $q^{-j}$.

\sm

The last sentence is Theorem 1.3.b. By Proposition \ref{pr:P-fixed},
any $\bGL^0$-invariant subspace in $L^2$ is generated by $\bfP$-fixed vectors. There
is only one $\bfP$-fixed vector in  each eigenspace of $\Delta$. 
Therefore
representations of $\bGL^0$ in eigenspaces of $\Delta$ are irreducible.

%%%%%%%%%%%%%%%%%%%%%%%%%%%%%%%%%%%%%%%%%%%%%%%%%%%%%%%%%
%%%%%%%%%%%%%%%%%%%%%%%%%%%%%%%%%%%%%%%%%%%%%%%%%%%%%%%%%%
%%%%%%%%%%%%%%%%%%%%%%%%%%%%%%%%%%%%%%%%%%%%%%%%%%%%%%%%%%

\section{Measures on flag manifolds}

\COUNTERS

Denote by $\bGr^\alpha$ the set of all subspaces of relative dimension $\alpha$.
The operator $J\in \bGL(\ell\oplus\ell^\diamond)$ send $\bGr^\alpha$ to $\bGr^{\alpha+1}$.

\sm

%%%%%%%%%%%%%%%%%%%%%%%%%%%%%%%%%%%%%%%%%%%%%%%%%%%%%%%%%

{\bf\punct Measures on finite flags.}
Let $L\in\bGr$.
Let the measures $\nu^\uparrow_L$ and $\nu^\downarrow_L$ be the same as above
(Lemma \ref{l:nu}). We regard them as measures on $\bGr$.
Now we define the measure $\nu^\uparrow_{[\alpha,\alpha+1]}$ on
$\bGr^\alpha\times\bGr^{\alpha+1}$ by
$$
\nu^\uparrow_{[\alpha,\alpha+1]}(S)=\int_{\bGr^\alpha} \nu^\uparrow_L\bigl(S\cap(\{L\}\times \bGr^{\alpha+1})\bigr)
\,d\mu(L)
$$ 
The measure $\nu^\uparrow_{[\alpha,\alpha+1]}$ is supported by the pairs of subspaces
\begin{equation}
L\in\bGr^\alpha,\quad K\in\bGr^{\alpha+1}, \qquad L\subset K.
\label{eq:two-step}
\end{equation}

We also define the measure 
$\nu^\downarrow_{[\alpha,\alpha-1]}$  on
$\bGr^{\alpha}\times\bGr^{\alpha-1}$ by
$$
\nu^\downarrow_{[\alpha,\alpha-1]}(S)=\int_{\bGr^{\alpha1}} \nu^\downarrow_L\bigl(S\cap(\{L\}\times \bGr^{\alpha})\bigr)
\,d\mu(L)
$$ 

\begin{lemma}
{\rm a)} The measure $\nu^\uparrow_{[\alpha,\alpha+1]}$ is a unique finite
$\bGL^0$-invariant measure supported by flags {\rm(\ref{eq:two-step})}.

\sm

{\rm b)} Transposition of factors in $\bGr^\alpha\times\bGr^{\alpha-1}$
send
$\nu^\downarrow_{[\alpha,\alpha-1]}$ to $\nu^\uparrow_{[\alpha-1,\alpha]}$.
\end{lemma}

{\sc Proof} is the same as proof of Lemma \ref{l:nu-symmetry}.
\hfill $\square$

\sm

For this reason we do not write arrows in superscripts
of $\nu_{[\alpha,\alpha+1]}$.

Next, we define a measure  $\nu_{[\alpha,\alpha+2]}$
on $\bGr^\alpha\times\bGr^{\alpha+1}\times\bGr^{\alpha+2}$
by
$$
\nu_{[\alpha,\alpha+2]}(S)=\int_{\bGr^\alpha\times\bGr^{\alpha+1}}
\nu^\uparrow_K\bigl(\{L\}\times\{K\}\times\bGr^{\alpha+2}\bigr)
\,d\nu_{[\alpha,\alpha+2]}(L,K)
.
$$
We get a measure supported by the space of flags
$$
L\subset K\subset N,\qquad L\in\bGr^\alpha, K\in\bGr^{\alpha+1}, M\in \bGr^{\alpha+2}
.$$

Iterating this operation we get measures $\nu_{[\alpha,\beta]}$
supported by the space $\bFl[\alpha,\beta]$ of flags
$$
L_\alpha \subset L_{\alpha-1}\subset\dots \subset L_\beta, \qquad L_j\in \bGr^j
$$
 By construction, this measure
is invariant with respect to the group $\bGL^0(\ell\oplus\ell^\diamond)$. 
Evidently it is the unique  invariant probabilistic measure. 

\sm

%%%%%%%%%%%%%%%%%%%%%%%%%%%%%%%%%%%%%%%%%%%%%%%

{\bf\punct Measures on a space of infinite flags.}
If a segment $[\gamma,\delta]$ contains a segment $[\alpha,\beta]$.
Then there is a forgetting map $\bFl[\gamma,\delta]\to\bFl[\alpha,\beta]$.
Since the  measure $\nu_{[\gamma,\delta]}$ is invariant, its pushforward is invariant
and therefore coincides with $\nu_{[\alpha,\beta]}$.
Therefore we get a $\bGL^0(\ell\oplus\ell^\diamond)$-invariant measure
on the space $\bFl(-\infty,\infty)$ of complete flags
$$
\dots \subset L_{-1}\subset L_0\subset L_1\subset\dots,
\qquad L_j\in\bGr^j.
$$

{\tt Math.Dept., University of Vienna,

 Nordbergstrasse, 15,
Vienna, Austria

\&

Institute for Theoretical and Experimental Physics,

Bolshaya Cheremushkinskaya, 25, Moscow 117259,
Russia

\&

Mech.Math.Dept., Moscow State University,

Vorob'evy Gory, Moscow

e-mail: neretin(at) mccme.ru

URL:www.mat.univie.ac.at/$\sim$neretin

wwwth.itep.ru/$\sim$neretin
}

\end{document}